\documentclass[a4paper]{article}

\usepackage{amsfonts}
\usepackage{bbm}
\usepackage{amssymb}
\usepackage{amsthm}
\usepackage{amsmath}
\usepackage{amscd}
\usepackage{graphicx}
\usepackage{stmaryrd}
\usepackage{scrpage2}
\usepackage[all]{xy}
\usepackage{tabularx}

\SelectTips{cm}{}

\DeclareMathOperator{\NIL}{NIL}
\DeclareMathOperator{\Nil}{Nil}
\DeclareMathOperator{\END}{END}
\DeclareMathOperator{\End}{End}
\DeclareMathOperator{\Mod}{Mod}
\DeclareMathOperator{\limit}{lim}
\DeclareMathOperator{\colim}{\ura{\limit}}
\DeclareMathOperator{\Coker}{Coker}
\DeclareMathOperator{\Ker}{Ker}
\DeclareMathOperator{\Image}{Im}
\DeclareMathOperator{\NK}{NK}
\newlength{\wlength}
\begin{document}

\title{The Behavior of Nil-Groups under Localization and the Relative Assembly Map}
\author{Joachim Grunewald\footnote{\sf email: grunewal@math.uni-muenster.de}}
\date{}

\newtheorem{Th}{Theorem}[section]
\newtheorem*{Th*}{Theorem}
\newtheorem{Pro}[Th]{Proposition}
\newtheorem*{Pro*}{Proposition}
\newtheorem{Le}[Th]{Lemma}
\newtheorem{Co}[Th]{Corollary}
\newtheorem*{Co*}{Corollary}

\theoremstyle{definition}

\newtheorem{De}[Th]{Definition}
\newtheorem*{De*}{Definition}
\newtheorem{Rem}[Th]{Remark}
\newtheorem{Ex}[Th]{Example}

\newcommand{\al}{\alpha}
\newcommand{\zz}{\mathbb{Z}}
\newcommand{\qq}{\mathbb{Q}}
\newcommand{\nn}{\mathbb{N}}
\newcommand{\ep}{\epsilon}
\newcommand{\Ga}{\Gamma}
\newcommand{\La}{\Lambda}
\newcommand{\la}{\lambda}
\newcommand{\xym}{\xymatrix}
\newcommand{\ra}{\rightarrow}
\newcommand{\Ra}{\Rightarrow}
\newcommand{\lra}{\longrightarrow}
\newcommand{\ot}{\otimes}
\newcommand{\op}{\oplus}
\newcommand{\ura}{\underrightarrow}
\newcommand{\tra}{\twoheadrightarrow}
\newcommand{\biggercup}{\ensuremath{\text{\Large$\bigcup$}}}
\newcommand{\tl}{\triangleleft}
\maketitle
\begin{center}
\mbox{ \textbf{Abstract}}
\end{center}
We study the behavior of the Nil-subgroups of $K$-groups
under localization. As a consequence of our results we obtain that the relative assembly map from the family of finite subgroups to the family of virtually cyclic subgroups is rationally an isomorphism. Combined with the equivariant Chern character we obtain a complete computation of the rationalized source of the $K$-theoretic assembly map that appears in the Farrell-Jones conjecture in terms of group homology and the $K$-groups of finite cyclic subgroups.
 
Specifically we prove that under mild assumptions we can always write the
Nil-groups and End-groups of the localized ring as a certain colimit over the
Nil-groups and End-groups of the ring, generalizing a result of Vorst. We define Frobenius and Verschiebung operations on certain Nil-groups. These operations provide the tool to prove
that Nil-groups are modules over the ring of Witt-vectors and are
either trivial or not finitely generated as abelian groups. Combining the localization results with the Witt-vector module structure we obtain that Nil and localization at an appropriate multiplicatively closed set $S$ commute, i.e. $S^{-1} \Nil = \Nil S^{-1}$. An important
corollary is that the Nil-groups appearing in
the decomposition of the $K$-groups of virtually cyclic groups are
torsion groups. 
\section{Introduction}
Over the last decades various kinds of Nil-groups appeared.
The most basic kind of Nil-group is given by Bass's definition of
the abelian groups $\Nil_i(R)$ \cite[Section 12.6]{Ba}. The so
called \emph{Fundamental Theorem of Algebraic $K$-Theory}
\cite{BHS,GQ} gives a description of the $K$-groups of the Laurent
polynomial ring of a ring $R$ in terms of the $K$-groups of $R$
and the groups~$\Nil_i(R)$ for all $i \in \zz$:
$$
K_i(R[t,t^{-1}]) = K_i(R) \op K_{i-1}(R) \op \Nil_{i-1}(R) \op
\Nil_{i-1}(R).
$$

If $\al$ is a ring automorphism of $R$, Farrell introduced the
abelian groups $\Nil_i(R;\al)$ \cite{Fa2}. He generalized,
together with Hsiang, the Fundamental Theorem of Algebraic
$K$-Theory to $K_1$ of the twisted Laurent polynomial ring
\cite{FH}. They proved that the sequence
$$
\xym{K_1(R) \ar[r]^-{1- \al_{*}} & K_1(R) \ar[r] &
K_1(R_{\al}[t,t^{-1}]) / (\Nil_0(R,\al) \op \Nil_0(R,\al^{-1}))
\ar[r] &}
$$
$$
\xym{ \ar[r] & K_0(R) \ar[r]^-{1- \al_{*}} & K_0(R) }
$$
is exact. This decomposition was extended to higher
algebraic $K$-theory by Grayson \cite{Gr1}.

Waldhausen introduced, for $R$-bimodules $X$ and $Y$, the abelian Nil-groups $\Nil_i(R;X,Y)$ \cite{Wa1,Wa4}. Nil-groups of this kind 
 appear in a long exact sequence relating the $K$-groups of a generalized free product to the $K$-groups of the ground rings. The sequence was extended to lower $K$-groups by Bartels and L\"uck \cite{BaL}. For $R$-bimodules $X$, $Y$, $Z$ and $W$ Waldhausen introduced the
abelian groups $\Nil_i(R;X,Y,Z,W)$ \cite{Wa1,Wa4}, which are the
most general kind of Nil-groups. Nil-groups of this kind appear in a long exact sequence relating 
the $K$-groups of a generalized Laurent extension to the
$K$-groups of the ground rings. Again, the sequence was extended to lower $K$-theory by Bartels and L\"uck \cite{BaL}. 
To avoid confusion, Nil-groups of the form $\Nil_i(R)$ are called \emph{Bass Nil-groups}, Nil-groups of the form $\Nil_i(R;\al)$ are called \emph{Farrell Nil-groups}, Nil-groups of the form~$\Nil_i(R;X,Y)$ are called \emph{Waldhausen Nil-groups of
generalized free products} and Nil-groups of the form
$\Nil_i(R;X,Y,Z,W)$ are called \emph{Waldhausen Nil-groups of
generalized Laurent extensions}.

All these Nil-groups have in common that they are defined as a subgroup of the $K$-theory of a certain Nil-category, which is in the following denoted by $\NIL(R)$, $\NIL(R;\al)$, $\NIL(R;X,Y)$ and $\NIL(R;X,Y,Z,W)$ respectively. 
\subsection*{The Behavior of Nil-Groups under Localization}
Nil-groups seem to be hard to compute. For example, it is known
that higher Bass Nil-groups are either trivial or not finitely
generated as abelian groups \cite{Fa,We1}. However, we prove that all kinds of
Nil-groups behave nicely under localization.
\begin{De*}
Let $R$ be a ring.
\begin{enumerate}
\item Let $T \subseteq R$ be a multiplicatively closed subset of
central non zero divisors. The ring $T^{-1} R$ is denoted by
$R_T$ and called the localization of $R$ at $T$.
\item  Let $s$ be a central non zero divisor and let $S$ be the
multiplicatively closed set generated by $s$. We use the short
hand notation $R_s$ for $R_S$.
\item Let $X$ be an $R$-bimodule. The
$R_T$-bimodule $R_T \ot_{R} X \ot_{R} R_T$ is denoted by ${}_T X_T$.
We use the short hand notation ${}_s X_s$ for $R_s \ot_R X \ot_R R_s$.
\end{enumerate}
\end{De*}

We prove the following result and similar results for the other kind of Nil-groups and End-groups:
\begin{Th*}
Let $R$ be a ring and let X, Y, Z and W be left flat
$R$-bimodules. Let $s$ be an element of the center of $R$ which is
not a zero divisor and satisfies $s \cdot x = x \cdot s$ for all
elements $x \in X$ and similar conditions for $Y$, $Z$ and $W$. We
obtain an isomorphism
$$ \zz [t,t^{-1}] \ot_{\zz[t]} \Nil_i(R;X,Y,Z,W) \cong
\Nil_i(R_s;{}_s X_s,{}_s Y_s,{}_s Z_s,{}_s W_s),$$
for $i \in \zz$, and $t$ acts on $\Nil_i(R;X,Y,Z,W)$ via the map
induced by the functor
\begin{align*}
\text{S} \colon \NIL(R;X,Y,Z,W) & \ra \NIL(R;X,Y,Z,W)\\
 (P,Q,p,q) & \mapsto (P,Q,p \cdot s,q \cdot s).
\end{align*}
Analogous results hold for the other kind of Nil-groups and for End-groups.
\end{Th*}

The condition that the bimodules $X$, $Y$, $Z$ and $W$ are left
flat does not seem to be overly restrictive since in all the cases
considered by Waldhausen $X$, $Y$, $Z$ and $W$ are left free by
the purity and freeness condition. The condition $s \cdot x = x
\cdot s$ translates in Waldhausen's setting of a generalized free product to the assumption that $s$ is mapped, under the
maps over which the pushout is formed, to central elements. The result was already known for Bass Nil-groups \cite{Vo}.
\subsection*{Operations on Nil-Groups}
For applications the result that localization and Nil commute is fruitful. To derive this result we develop Frobenius and Verschiebung operations on certain kinds of Nil-groups. For Bass Nil-groups and
$\End$-groups Frobenius and Verschiebung operations are well understood
\cite{Bl2,CS,Sti,We1}. 

Our applications of these operations are twofold. Firstly we generalize a non-finiteness result of Farrell \cite{Fa}.
\begin{Co*}
Let $R$  be a ring, let $G$ be a group, let $X$ and $Y$ be
arbitrary $RG$-bimodules and let $\al$ and $\beta$ be inner group
automorphisms of $G$. Then $\Nil_i(RG; \al)$ and
$\Nil_i(RG;RG_{\al} \op X ,RG_{\beta} \op Y)$ are either trivial
or not finitely generated as an abelian group for $i \in \zz$.
\end{Co*}

Secondly they are the main tool for the proof that Nil-groups are modules over the ring of Witt vectors.
\begin{Pro*}
Let $R$ be a commutative ring, $G$ a group and $\al,
\beta \colon G \ra G$  inner group automorphism. The groups
$\Nil_i(RG;\al)$ and~$\Nil_i(RG;RG_{\al},RG_{\beta})$ are continuous modules
over the ring of Witt vectors of $R$.

Moreover the Witt vector-module structure is compatible with the Frobenius and Verschiebung operations in the following sense: If for a natural number $n$ the Frobenius operation is denoted by $\text{F}_n$, the Verschiebung by $\text{V}_n$ and the Witt vector-module structure by $\ast$ we have
\[
\text{V}_n(y \ast \text{F}_n x) = (V_n y) \ast x
\]
for every element $y$ of the ring of Witt vectors and $x$ an element in the Nil-group.
\end{Pro*}

To get this Witt vector-module structure on Nil-groups, we first
define an $\End_0(R)$-module structure on Nil-groups. It is a result of Almkvist that $\End_0(R)$ is a dense subring of the ring of Witt vectors \cite{Al2}. To prove that the $\End_0(R)$-module structure can be
extended to a Witt vector-module structure we use the Frobenius and Verschiebung operations. 

Combined with the nice behavior of the Nil-groups under localization the Witt vector-module structure implies that Nil and localization commutes:
\begin{Th*}
Let $R$ be $\zz_T$ for some multiplicatively closed set $T
\subseteq \zz -\{ 0 \}$, $\hat{\zz}_p$ or a commutative
$\qq$-algebra. Let $G$ be a group and let $\al$ and $\beta$ be
inner automorphisms of $G$. Then for every multiplicatively closed
set $S \subset R$ of non zero divisors there are isomorphisms of
$R_S$-modules
$$
 R_S \ot_{R} \Nil_i(RG;\al) \cong \Nil_i(R_S G;\al)
$$
and
$$
 R_S \ot_{R} \Nil_i(RG;RG_{\al},RG_{\beta})  \cong
 \Nil_i(R_S G;R_S G_{\al}, R_S G_{\beta}),
$$
for all $i \in \zz$.
\end{Th*}
Note that we do not need to assume $S$ to be central since all considered rings are commutative. An important part of the paper at hand derives this theorem from the theorem given above by using Witt vector techniques. Observe that the second theorem is not an immediate corollary of the first theorem.
\subsection*{Torsion Results}
The preceding theorem implies, combined with induction and transfer maps on the Nil-groups, the following torsion results.  

Every polycyclic-by-finite group $G$ contains a poly-infinite cyclic subgroup of finite index.
\begin{Th*}
Let G be a polycyclic-by-finite group containing a poly-infinite cyclic subgroup of finite index $n$. Let $\al$, $\beta$ and $\gamma$ be group automorphisms such that $\al$ is of finite order $m$ and $\beta \circ \gamma$ is of finite order $m'$.
\begin{enumerate}
\item The group $\Nil_i(\zz G;\al)$ is an $(n \cdot m)$-primary torsion
group for $i \in \zz$.
\item The group $\Nil_i(\zz G;\zz G_{\beta},\zz G_{\gamma})$ is an
$(n \cdot m')$-primary torsion group for $i \in \zz$.
\end{enumerate}
\end{Th*}

As an important application of this theorem, we get that the
Nil-groups appearing in the calculation of the $K$-groups of infinite
virtually cyclic groups are torsion groups. There are two kinds of infinite virtually cyclic groups:
\begin{enumerate}
\renewcommand{\theenumi}{\Roman{enumi}}
\renewcommand{\labelenumi}{\theenumi.}
\item the semidirect product $G \rtimes \zz$ of a finite group $G$
and the infinite cyclic group;
\item the amalgamated product $G_1 \ast_H G_2$ of two finite
groups $G_1$ and $G_2$ over a subgroup~$H$ such that $[G_1:H] =
2=[G_2:H]$.
\end{enumerate}

In the case of a virtually cyclic group of the first type
Farrell Nil-groups of finite groups appear. If we consider a
virtually cyclic group of the second type the Nil-groups
$\Nil_i(\zz H,\zz [G_1 - H],\zz [G_2 - H])$ relate the $K$-groups
of $G_1 \ast_H G_2$ to the $K$-groups of $H$, $G_1$ and $G_2$. The
group $H$ is an index two subgroup. Thus we can find automorphisms
$\al$ and $\beta$ of $H$ such that the $\zz H$-bimodules $\zz [G_1 -
H]$ and $\zz [G_2 - H]$ are isomorphic to $\zz H_{\al}$ and $\zz
H_{\beta}$. 
\begin{Co*}
Let G be a finite group of order n and let $\al$ and $\beta$ be
group automorphisms such that $\al \circ \beta$ is of finite order $m$. The group
$\Nil_i(\zz G;\zz G_{\al},\zz G_{\beta})$ is an $(n \cdot m)$-primary torsion group for $i \in \zz$.
\end{Co*}

The result that Nil-groups of finite groups are torsion groups was
already known in other cases: For Bass Nil-groups Weibel proved
that if $G$ is a finite group of order~$n$, then $\Nil_i(\zz G)$
is $n$-primary torsion for $i \geq 0$ \cite{We1}. Our approach is similar to his. If $i \leq -1$, Bass Nil-groups of $\zz G$ are known to vanish \cite{Ba,We2}. Farrell and Jones proved that the groups $\Nil_{i}(\zz G;\al)$ and $\Nil_i(\zz G;\zz G_\al,\zz G_\beta)$ are trivial for $i \leq -2$ \cite{FJ2}. Connolly and Prassidis took a similar approach to prove that
$\Nil_{-1}(\zz G;\al)$ and $\Nil_{-1}(\zz G; \zz G_\al, \zz G_\beta)$ are torsion \cite{CP}. Kuku and Tang generalized this
concept to prove that $\Nil_i(\zz G,\al)$ is $n$-primary torsion 
for $i \geq -1$ and that $\Nil_0(\zz H,\zz [G_1 - H] ,\zz [G_2 - H])$ is $n$-primary torsion \cite{KT}. Note that in Kuku and Tang's paper the Nil-groups appearing in the decomposition of Waldhausen are denoted by
$\widetilde{\Nil}_{i}^W\negthinspace(R;R^{\al},R^{\beta})$ and
called Waldhausen Nil-groups. It was not known that $\Nil_i(\zz H,\zz [G_1 - H] ,\zz [G_2 - H])$ are torsion groups for $i$ bigger than zero, as was
incorrectly stated in~\cite{LR}.

For an arbitrary group $G$ it is unknown whether the groups $\Nil_i(\qq G; \al)$ and $\Nil_i(\qq G;\qq G_{\al},\qq G_{\beta})$ are trivial. Another important application of our results is that 
these groups are almost torsion free:
\begin{Th*}
Let G be an arbitrary group and let $\al$, $\beta$ and $\gamma$ be group
automorphisms such that $\al$ is of finite order $m$ and $\beta \circ \gamma$ is of finite order $m'$.
\begin{enumerate}
\item The group $\zz[1/m] \ot_{\zz} \Nil_i(\qq G;\al)$ is a
$\qq$-module for $i \in \zz$. 
\item The group $\zz[1/m'] \ot_{\zz} \Nil_i(\qq G;\qq
G_{\beta},\qq G_{\gamma})$ is a $\qq$-module for $i \in \zz$. 
\end{enumerate}
\end{Th*}
\subsection*{The Farrell-Jones Conjecture}
The Farrell-Jones conjecture predicts that the \emph{assembly map} 
\[
A_{\mathcal{VC}yc \ra \mathcal{A}ll} \colon
H_i^G(E_{\mathcal{VC}cy}(G); \textbf{K}_{R} ) \ra  H_i^G (
E_{\mathcal{A}ll}(G) ; \mathbf{K}_R) \cong K_i(R G)
\]
is an isomorphism. For a survey on the Farrell-Jones conjecture see for example~\cite{LR}.

We can study the left hand side of the Farrell-Jones conjecture by subfamilies $\mathcal{F}$ of the family of virtually cyclic groups. The smaller the
family $\mathcal{F}$ is, the easier it is to compute
$H_i^G(E_{\mathcal{F}}(G); \mathbf{K}_{R})$. It is known that the \emph{relative assembly map} from the family of finite subgroups $\mathcal{F}in$ to the family of virtually cyclic subgroups $\mathcal{VC}yc$
\[
A_{\mathcal{F}in \ra \mathcal{VC}yc} \colon H_i^G(
E_{\mathcal{F}in}(G) ; \mathbf{K}_R ) \ra H_i^G( E_{\mathcal{VC}yc}(G);
\mathbf{K}_R)
\]
is split injective \cite{Bar}. In Theorem \ref{mainth}, we prove that
rationally this relative assembly map is an isomorphism. The main
ingredient of the proof is that Nil-groups of finite groups are
torsion. Combining this result with L\"{u}ck's \cite{Lue2} computation of $H_i^G(E_{\mathcal{F}in}(G);\textbf{K}_R ) \ot \qq $ we obtain a complete computation of the rationalized source of the assembly map in terms of group homology and the $K$-theory of finite cyclic subgroups. Before we state the result let us recall the relevant notions. For a finite group $G$ the groups $K_i(R G)$ are modules over the Burnside ring $A(G)$. We have a ring homomorphism 
\begin{align*}
\chi^G \colon A(G) & \ra \prod_{(H) \in \mathcal{F}in} \zz \\
[S] & \mapsto (|S^H|)_{(H) \in \mathcal{F}in}
\end{align*}
which sends the class of a finite $G$-set $S$ to the element given by the cardinality of the $H$-fixed point set. The map $\chi^G$ becomes rationally an isomorphism
\[
\chi^G \ot \text{id} \colon A(G) \ot_\zz \qq \ra \prod_{(H) \in \mathcal{F}in} \qq
\]
and we define $\theta_G \in A(G) \ot_\zz \qq$ to be the element that is mapped under $\chi^G \ot \text{id}$ to the element, whose entry is one if $(H)=(G)$ and zero otherwise.  
\begin{Co*} For every group $G$, the source of the rationalized assembly map is
\[
H_i^G(E_{\mathcal{VC}yc}(G) ; \mathbf{K}_\zz ) \ot \qq \cong \bigoplus_{p+q=i} \bigoplus_{(C) \in (\mathcal{FC}y)} H_p(B Z_G C;\qq) \ot_{\qq[W_G C]} \theta_C \cdot K_q(\zz C) \ot \qq
\]
where $(\mathcal{FC}y)$ is the set of conjugacy classes $(C)$ of finite cyclic subgroups, $N_G H$ is the normalizer of a subgroup $H$, $Z_G H$ is the centralizer and $W_G H := N_G H / (H \cdot Z_GH)$. 
\end{Co*}
The Farrell-Jones conjecture is known to be true for a large class of groups. For all groups $G$ for which the Farrell-Jones conjecture is known to be true the preceding corollary implies:
\[
K_i(\zz G) \ot \qq \cong \bigoplus_{p+q=i} \bigoplus_{(C) \in (\mathcal{FC}y)} H_p(B Z_G C;\qq) \ot_{\qq[W_G C]} \theta_C \cdot K_q(\zz C) \ot \qq. 
\]
\subsection*{Acknowledgments}
This paper is essentially the author's Ph.D. thesis, which was written under the guidance of Wolfgang L\"{u}ck and Holger Reich at the University of M\"{u}nster \cite{Phdthe}. It is a pleasure to thank them for countless valuable discussions. Finally, I would like to thank the referee for his careful reading of the manuscript. 
\section{End and Nil-groups}
Various kinds of End- and Nil-groups have been defined \cite{Ba,Fa2,Gr4,Wa1,Wa4} as subgroups of the $K$-groups of
End- and Nil-categories. In the following section, we give a unified definition of End- and Nil-groups and of End-
and Nil-categories. All Nil-groups from the introduction will appear as spacial cases.
\subsection{End-Groups and End-categories}
We begin by generalizing End-groups and End-categories.
\begin{De}[End-category] \label{End}
Let $\mathcal{A}$ be an abelian category, let $\text{F} \colon
\mathcal{A} \ra \mathcal{A}$ be an exact functor and let $\mathcal{C}
\subseteq \mathcal{A}$ be a full subcategory which is \emph{closed
under extension}, i.e., if
$$
\xym{0 \ar[r] & A_1 \ar[r] & A_2 \ar[r] & A_3 \ar[r] & 0}
$$
is exact and both $A_1$ and $A_3$ belong to $\mathcal{C}$ then
$A_2$ belongs to $\mathcal{C}$. We define $\END(\mathcal{C};\text{F})$ to
be the following category. Objects are pairs $(C,c)$ consisting of
an object $C$ of $\mathcal{C}$ and a morphism
$$
c \colon C \ra \text{F}(C)
$$
in $\mathcal{A}$. A morphism $f \colon (C,c) \ra (C',c')$ in
$\END(\mathcal{C};\text{F})$ consist of a morphism $f \colon C \ra
C'$ in $\mathcal{C}$ satisfying $c' \circ f = \text{F}(f) \circ
c$.

Let $\text{S} \colon \END(\mathcal{C}; \text{F}) \ra \mathcal{C}$
be the forgetful functor mapping $(C,c)$ to $C$ and a morphism $f
\colon (C,c) \ra (C',c')$ to $f \colon C \ra C'$. A sequence
$$
\xym{ (C_1,c_1) \ar[r]^{f} & (C_2,c_2) \ar[r]^{g} & (C_3,c_3) }
$$
is called \emph{exact} at $(C_2,c_2)$ if the sequence
$$
\xym{ C_1 \ar[r]^{f} & C_2 \ar[r]^{g} & C_3 }
$$
in $\mathcal{C}$ obtained by applying S is exact at $C_2$. Short
exact sequences, surjective and injective morphism in
$\END(\mathcal{C};\text{F})$ are defined in the obvious way.
\end{De}
\begin{Pro}
Let the notation be as in the preceding definition. We
additionally assume that $\mathcal{C}$ has a small skeleton. The
category $\END(\mathcal{A};\text{F})$ is an abelian category and
the category $\END(\mathcal{C};\text{F})$ is an exact category.
\end{Pro}
\begin{proof}
One easily verifies the required identities.
\end{proof}
\begin{De}[End-groups]
Let the notation be as in the preceding definition.
We define
$$
\End_i(\mathcal{C};\text{F}):= \Ker \big( K_i(\text{S}) \colon
K_i\big(\END(\mathcal{C};\text{F})\big) \ra K_i(\mathcal{C})
\big)
$$
for $i \geq 0$.
\end{De}
Consider the functor 
\begin{align*}
\text{T} \colon \mathcal{C} & \rightarrow \END(\mathcal{C};\text{F}) \\
C & \mapsto (C,0).
\end{align*}
Since $\text{S} \circ \text{T} = \text{Id}$ we have
\[
 K_i(\END(\mathcal{C};\text{F})) = K_i(\mathcal{C}) \op \End_i(\mathcal{C};\text{F}).
\]

We continue by defining maps between End-groups. Suppose that we
have for $j = 0, 1$ abelian categories $\mathcal{A}_j$ together
with full subcategories $\mathcal{C}_j \subseteq \mathcal{A}_j$
that are closed under extension. Additionally we have exact functors $\text{F}_j
\colon \mathcal{A}_j \ra \mathcal{A}_j$. Suppose that $u \colon
\mathcal{A}_0 \ra \mathcal{A}_1$ is a functor that sends $\mathcal{C}_1$ to $\mathcal{C}_2$ and $U \colon u
\circ \text{F}_0 \ra \text{F}_1 \circ u$ is a natural
transformation of exact functors. Then we obtain an
exact functor
$$
\END(u,U) \colon \END(\mathcal{C}_0; \text{F}_0) \ra
\END(\mathcal{C}_1; \text{F}_1)
$$
which sends an object $c \colon C \ra \text{F}_0(C)$ to the object
given by the composition
$$
\xym{u(C) \ar[r]^-{u(c)} & u(\text{F}_0(C)) \ar[r]^{U(C)} &
\text{F}_1(u(C)).}
$$
A morphism $f \colon (C,c) \ra (C',c')$ is sent to the morphism
whose underlying morphism in $\mathcal{C}$
is $u(f) \colon u(C) \ra u(C')$. This is a well-defined functor
since the following diagram commutes.
$$
\xym{u(C) \ar[r]^-{u(c)} \ar[d]^{u(f)} & u(\text{F}_0(C))
\ar[r]^{U(C)} \ar[d]^{u(\text{F}_0(f))} & \text{F}_1(u(C))
\ar[d]^{\text{F}_1(u(f))} \\
u(C') \ar[r]^-{u(c')} & u(\text{F}_0(C')) \ar[r]^-{U(C')} &
\text{F}_1 (u(C')).}
$$
Since also the diagram
$$
\xym@C=3.5em{\END(\mathcal{C}_0;\text{F}_0) \ar[r]^{\END(u,U)}
\ar[d]_{\text{S}} &
\END(\mathcal{C}_1;\text{F}_1) \ar[d]^{\text{S}} \\
\mathcal{C}_0 \ar[r]^u & \mathcal{C}_1}
$$
commutes, we obtain for the given pair $(u,U)$ a homomorphism
$$
\End_i(u,U) \colon \End_i(\mathcal{C}_0;\text{F}_0) \ra
\End_i(\mathcal{C}_1;\text{F}_1).
$$

Let $R$ be a ring. We define $\Mod (R)$ to be the abelian category of all
right modules and $\textbf{P}(R)$ to be the exact subcategory of all
finitely generated projective right $R$-modules.
We denote the identity functor on $\Mod (R)$ by Id. 
\begin{De}[Higher End-groups]
Let $R$ be a ring. For $i \geq 0$, we define
$$
\End_i(R) := \End_i(\mathbf{P}(R);\text{Id}).
$$
\end{De}

The \emph{cone ring} $\La \zz$ of $\zz$ is the ring of matrices
over $\zz$ such that every column and every row contains only
finitely many non zero entries. The \emph{suspension ring} $\Sigma
\zz$ is the quotient of $\La \zz$ by the ideal of finite matrices.
For a natural number $n \geq 2$ we define inductively $\Sigma^n
\zz= \Sigma \Sigma^{n-1} \zz$. For an arbitrary ring $R$ we define
the $n$-fold suspension ring $\Sigma^n R$ to be $\Sigma^n \zz \ot
R$. For an $R$ bimodule $X$ we define $\Sigma^n X$ to be the
$\Sigma^n R$-bimodule $\Sigma^n R \ot X$ where the right $\Sigma^n
R$-module structure is given by $(z_1 \ot x) \cdot (z_2 \ot r) = z_2 \cdot z_1 \ot  x \cdot r$ for $x \in X$, $z_1$, $z_2 \in \Sigma^n \zz$ and $r \in R$.
\begin{De}[Lower End-groups]
Let $R$ be a ring. For $i < 0$, we define
$$
\End_i(R) := \End_0(\mathbf{P}(\Sigma^{-i} R);\text{Id}).
$$
\end{De}
\subsection{Nil-Groups and Nil-categories}
In the sequel the full subcategory of $\END(\mathcal{C};\text{F})$
consisting of nilpotent endomorphisms will become important.
\begin{De}[Nil-category]\label{NilCategory}
Let the notation be as in Definition \ref{End}, let $(C,c)$ be an object of $\END (\mathcal{C};\text{F})$ and for a natural number $n \geq 0$ let F$^n$ be the $n$-fold composite $\text{F} \circ \text{F} \cdots \circ \text{F}$.
We define $c^n$ to be the morphism given by the composite
$$
\xym@C=4em{ \text{F}^0(C)=C \ar[r]^-{\text{F}^0(c)=c} &
\text{F}(C) \ar[r]^-{\text{F}(c)} & \text{F}^2(C)
\ar[r]^-{\text{F}^2(c)} & \cdots \ar[r]^-{\text{F}^{n - 1}(c)} &
\text{F}^n (C). }
$$
The object $(C,c)$ is called \emph{nilpotent} if for
some natural number $N \geq 1$ the morphism $c^N$
is given by the zero morphism from $C$ to $\text{F}^N (C)$. We
call the smallest number with this property the \emph{nilpotency
degree} of $(C,c)$. Let $\NIL(\mathcal{C};\text{F})$ be the full
subcategory of $\END(\mathcal{C};\text{F})$ given by nilpotent
objects.
\end{De}
\begin{Pro} \label{5}
Let the notation be as in the preceding definition. We
additionally assume that $\mathcal{C}$ has a small skeleton. The
category $\NIL(\mathcal{A};\text{F})$ is an abelian category and
the category $\NIL(\mathcal{C};\text{F})$ is an exact category.
\end{Pro}
\begin{proof}
One easily verifies the required identities.
\end{proof}
\begin{De}[Nil-groups]
Let the definition be as in Definition \ref{End}.
The functor S defined above restricts to a functor on $\NIL (\mathcal{C};\text{F})$. We define
$$\Nil_i(\mathcal{C};\text{F}):= \Ker \big( K_i(\text{S}) \colon
K_i\big(\NIL(\mathcal{C};\text{F})\big) \ra K_i (\mathcal{C})
\big)$$
for $i \geq 0$.
\end{De} 
Similarly as above we have
\[
  K_i\big(\NIL(\mathcal{C};\text{F})\big) = K_i(\mathcal{C}) \op  \Nil_i(\mathcal{C};\text{F})
\]
and we obtain for a pair $(u,U)$ a morphism
$$
\Nil_i(u,U) \colon \Nil_i(\mathcal{C}_0;\text{F}_0) \ra
\Nil_i(\mathcal{C}_1;\text{F}_1).
$$

In the following we will show how the Nil-groups from the introduction fit into the given setting. The definition of lower Nil-groups using the suspension ring is one possible approach. One could also use results of Carter \cite{Car} or Schlichting \cite{Schl}.
\begin{De}[Bass Nil-groups]
Let $R$ be a ring. We define \emph{Bass Nil-groups}, for $i \geq
0$, by
$$
\Nil_i(R) := \Nil_i(\mathbf{P}(R);\text{Id})
$$
and for $i < 0$ by
$$
\Nil_i(R) := \Nil_0(\mathbf{P}(\Sigma^{-i} R);\text{Id}).
$$
\end{De}
\begin{Rem}
Notice that the given definition for lower Nil-groups coincides
with the definition of lower $\NK$-groups given by Bass \cite{Ba}. This follows from the fact that for an arbitrary ring $R$ and a natural number $i < 0$
we have $K_i(R) = K_0(\Sigma^{-i} R)$ and $\Sigma R [t] = (\Sigma
R )[t]$.
\end{Rem}
\begin{De}[Farrell Nil-groups]
Let $R$ be a ring, $X$ a left flat $R$-bimodule and $\text{F}_X$ the exact functor from $\Mod (R)$ to $\Mod (R)$ which is induced by tensoring with $X$ on the right.

We define \emph{Farrell Nil-groups}, for $i \geq 0$, by
$$
\Nil_i(R;X) := \Nil_i(\mathbf{P}(R);\text{F}_{X})
$$
and for $i < 0$ by
$$
\Nil_i(R;X) := \Nil_0(\mathbf{P}(\Sigma^{-i} R);\text{F}_{\Sigma^{-i} X}).
$$
\end{De}
\begin{Rem}\label{Rem}
\begin{enumerate}
\item  The bimodules appearing in the
decomposition of Farrell and Hsiang are of the form $X=R$ where
the left $R$-module structure is given by multiplication and the
right $R$-module structure comes from an automorphism $\al$ of the
ring $R$. These kind of bimodules are denoted by $R_\al$.
\item Notice that the given definition of lower Farrell Nil-groups
coincides with the usual definition of lower Farrell Nil-groups.
Again, this follows from the fact that for an arbitrary ring $R$ and an
arbitrary endomorphism $\al$ of $R$ we have $\Sigma R_{\al} [t] =
(\Sigma R )_{\text{id} \ot \al}[t]$.
\end{enumerate}
\end{Rem}
\begin{De}[Waldhausen Nil-groups of generalized free products]
\label{hochdef} Let $R$ be a ring and let $X$ and $Y$ be left flat
$R$-bimodules. Define
$$
\text{F}_{X,Y} \colon \Mod (R) \times \Mod (R) \ra
\Mod (R) \times \Mod (R)
$$
by sending $(M,N)$ to $(N \ot X, M \ot Y)$. 

We define \emph{Waldhausen Nil-groups of generalized free
products}, for $i \geq 0$, by
$$
\Nil_i(R;X,Y): = \Nil_i(\mathbf{P}(R) \times \mathbf{P}(R);\text{F}_{X,Y})
$$
and for $i < 0$ by
$$
\Nil_i(R;X,Y): = \Nil_0(\mathbf{P}(\Sigma^{-i} R) \times \mathbf{P}(\Sigma^{-i} R);\text{F}_{\Sigma^{-i} X, \Sigma^{-i} Y}).
$$
\end{De}
\begin{De}[Waldhausen Nil-groups of generalized Laurent extension]
Let $R$ be a ring and let $X$, $Y$, $Z$ and $W$ be left flat
$R$-bimodules. Define
$$
\text{F}_{X,Y,Z,W} \colon \Mod (R) \times \Mod (R) \ra
\Mod (R) \times \Mod (R)
$$
by sending $(M,N)$ to $(N \ot X \op M \ot Z, M \ot Y \op N \ot
W)$. We define \emph{Waldhausen Nil-groups of generalized Laurent extensions},
for $i \geq 0$, by
$$
\Nil_i(R;X,Y,Z,W): = \Nil_i(\mathbf{P}(R) \times \mathbf{P}(R);\text{F}_{X,Y,Z,W})
$$
and for $i < 0$ by
$$
\Nil_i(R;X,Y,Z,W): = \Nil_0(\mathbf{P}(\Sigma^{-i} R) \times \mathbf{P}(\Sigma^{-i} R) ; \text{F}_{\Sigma^{-i} X,\Sigma^{-i} Y,\Sigma^{-i} Z, \Sigma^{-i} W}).
$$
\end{De}

In the sequel we will treat all Nil-groups at once. For this purpose we define F to be either of the functors Id, F$_X$, F$_{X,Y}$ or F$_{X,Y,Z,W}$. By abuse of languages we will denote the Nil-groups and categories by $\Nil(\mathbf{P}(R);\text{F})$ even thought for F$=\text{F}_{X,Y}$ or $=\text{F}_{X,Y,Z,W}$ it should be denoted by $\Nil(\mathbf{P}(R) \times \mathbf{P}(R);\text{F})$.
\section{The Behavior of Nil-Groups under Localization}
In this section we study the behavior of Nil-groups under localization. First we provide homological facts of $\NIL$-categories. Secondly we develop a long exact sequence for certain $K$-groups. Thirdly we use these results to develop a long exact localization sequence for the $K$-theory of $\NIL$-categories which is similar to the long exact localization sequence of algebraic $K$-theory \cite{GQ}. Our proof of the exactness follows an approach given by Grayson \cite{Gr3} and Staffledt \cite{St}. We use this sequence to obtain the localization results stated in the introduction.
\subsection{Homological Algebra of $\NIL$-categories}\label{Lifts}
\begin{Le}
\label{12} Let $\mathcal{A}$ be an abelian category, let
$\mathcal{C}$ be a subcategory which is closed under extension and let $\text{F} \colon \mathcal{A} \ra \mathcal{A}$ be an exact
functor. The subcategory $\NIL(\mathcal{C}; \text{F})$ of
$\END(\mathcal{C}; \text{F})$ is closed under extension.
\end{Le}
\begin{proof}
The results follows by a diagram chase.
\end{proof}
\begin{Le}\label{7} Let the notation be as in the preceding lemma.
Let
$$
 \xym@C=2.7em{0 \ar[r] & (C_1,c_1) \ar[r] & (C_2,c_2)
 \ar[r] & (C_3,c_3) \ar[r] & 0}
$$
be a short exact sequence in $\NIL (\mathcal{C};\text{F})$ and let
$(P_1,p_1)$ and $(P_2,p_2)$ be objects in
$\NIL(\mathcal{C};\text{F})$ admitting surjections
\begin{align*}
\pi_{P_1} \colon (P_1,p_1) & \twoheadrightarrow (C_1,c_1) \\
\pi_{P_2} \colon (P_2,p_2) & \twoheadrightarrow (C_3,c_3).
\end{align*}
If $P_1$ and $ P_2$ are projective, we can construct a surjection
\begin{align*}
\pi_P \colon P_1 \op P_2 & \twoheadrightarrow C
\end{align*}
out of $\pi_{P_1}$ and $\pi_{P_2}$.

There is an object $(P_1 \op P_2, p)$ in
$\NIL(\mathcal{C};\text{F})$ such that the diagram
$$
\xym{ 0 \ar[r] & (P_1,p_1) \ar[r] \ar@{->>}[d]^-{\pi_{P_1}} & (P_1 \op
P_2, p) \ar[r] \ar@{->>}[d]^{\pi_{P}}
& (P_2,p_2) \ar[r] \ar@{->>}[d]^{\pi_{P_2}} & 0 \\
0 \ar[r] & (C_1,c_1) \ar[r] & (C_2,c_2)
 \ar[r] & (C_3,c_3) \ar[r] & 0
 }
$$
commutes.
\end{Le}
\begin{proof}
Use the projectivity of $P_2$ to construct a map $p \colon P_1 \op P_2 \ra \text{F} (P_1 \op P_2)$ which fits into the given diagram. The map $p$ is nilpotent by the preceding lemma.
\end{proof}
\begin{Le}
Assume the following conditions for an abelian category
$\mathcal{A}$ with subcategory $\mathcal{C}$ which is closed under
extension and an exact functor $\text{F}$:
\begin{enumerate}
\item For any object $C$ in $\mathcal{C}$ there exists an object
$P$ in $\mathcal{C}$ which is projective in $\mathcal{A}$ and admits an epimorphism $f \colon P \ra C$;
\item Any object in $\mathcal{A}$ can be written as a colimit of a
directed system $\{ C_i | i\in I \}$ such that each structure map
is a monomorphism and each object belongs to $\mathcal{C}$;
\item The subcategory $\mathcal{C}$ is closed under quotient objects in $\mathcal{A}$, i.e. for any epimorphism $g \colon C \ra C'$ in $\mathcal{A}$ for
which $C$ belongs to $\mathcal{C}$ also $C'$ belongs to
$\mathcal{C}$;
\item Suppose that an object $A$ in $\mathcal{A}$ is the colimit
of a directed system $\{C_i | i \in I\}$ such that each structure
map is a monomorphism and each object belongs to $\mathcal{C}$. Let $h
\colon C \ra A$ be an injective morphism with $C \in
\mathcal{C}$. Then there exists an index $i \in I$ such that the
image of $h$ is contained in the image of $C_i \ra A$;
\item The functor $\text{F}$ commutes with colimits over directed
systems and structure maps which are monomorphisms.
\end{enumerate}
Then we can find for any object $(C,c)$ in
$\NIL(\mathcal{C};\text{F})$ an object $(P,p)$ in
$\NIL(\mathcal{C};\text{F})$ together with an epimorphism from
$(P,p)$ onto $(C,c)$ such that $P$ is projective. The nilpotency
degree of $(P,p)$ is smaller or equal to one plus the nilpotency
degree of $(C,c)$.
\end{Le}
\begin{proof}
The result is proven by induction on the nilpotency degree $L$.
For $L=1$ we have $c=0$. Let $P_C$ be a projective object
surjecting onto $C$. The trivial map from $P_C$ to $\text{F}(P_C)$
is a lift of $c$.

For the induction step from $L$ to $L + 1$ let $(C,c)$ be an
object of $\NIL (\mathcal{C};\text{F})$ of nilpotency degree $L +
1$. Let $K$ be the kernel of $c^{\ell}$. Since F is exact, the following
sequence is also exact.
$$
\xym{0 \ar[r] & \text{F}(K) \ar[r]^-{\text{F}(\iota)} & \text{F}(C)
\ar[r]^-{\text{F}(c^{\ell})} & \text{F}^{\ell+1}(C) \ar[r] & 0.}
$$
The image of $c$ is contained in the kernel of $c^{\ell}$, since
$c^{\ell +1}=0$. The given exact sequence implies that the image
of $c$ is contained in the image of $\text{F}(\iota) \colon
\text{F}(K) \ra \text{F}(C)$. By assumption we can write
$$
K = \colim_{i \in I} K_i
$$
for a directed system with injective structure maps such that each
$K_i$ belongs to $\mathcal{C}$. By assumption the canonical map
$$
 \colim_{i \in I} \text{F}(K_i) = \text{F}(K)
$$
is bijective. The image of $c \colon C \ra \text{F}(C)$ belongs,
by assumption, to $\mathcal{C}$ since $C$ belongs to
$\mathcal{C}$. Thus we can find an index $i_0 \in I$ such that the
image of $c$ is a subobject of the image of $\text{F}(K_{i_0}) \ra
\text{F}(K)$. Put $C_{\Image} = K_{i_0}$.

The object $C_{\Image}$ has the property that the image of $c$ is
contained in $\text{F}(C_{\Image})$. This has two implications.
First of all we can restrict $c$ to a morphism of $C_{\Image}$.
The object $(C_{\Image},c|_{C_{\Image}})$ is of nilpotency degree
$L$. By construction, $C_{\Image}$ is an object in
$\NIL(\mathcal{C};\text{F})$. Thus, by our induction hypothesis,
we get a projective object $P_{C_{\Image}}$ and a nilpotent lift
$\tilde{c}_{\Image}$ such that the diagram
$$
\xym{P_{C_{\Image}} \ar[r]^-{\tilde{c}_{\Image}} \ar@{->>}[d] &
\text{F} (P_{C_{\Image}}) \ar@{->>}[d]
\\ C_{\Image} \ar[r]^-{c|_{C_{\Image}}} \ar@{^{(}->}[d] &
\text{F} (C_{\Image})\ar@{^{(}->}[d] \\
C \ar[r]^-c & \text{F}(C)}
$$
commutes.

Let $P_C$ be a projective object surjecting onto $C$. The second
implication is that since $P_C$ is projective there is a map
$\tilde{c}_C$ making the diagram
$$
\xym{P_C \ar[r]^-{\tilde{c}_C} \ar@{->>}[d] &
\text{F}(P_{C_{\Image}}) \ar@{->>}[d]
\\ C \ar[r]^-c & \text{F}(C_{\Image})}
$$
commutative. Plugging these maps together in the matrix
$$
\tilde{c}:= \left(
\begin{array}{cc}
 0 & 0 \\
 \tilde{c}_C & \tilde{c}_{\Image}
 \end{array}\right)
$$
we get a commutative diagram
$$
\xym{P_C \op P_{C_{\Image}} \ar[r]^-{\tilde{c}} \ar@{->>}[d] &
\text{F}(P_C \op P_{C_{\Image}}) \ar@{->>}[d] \\
C \ar[r]^-c & \text{F}(C).}
$$

This lift of $c$ is nilpotent since $\tilde{c}_{\Image}$ is. Thus
the object $\big(P_C \op P_{C_{\Image}},\tilde{c}\big)$ is a
projective object of nilpotency degree $L+2$ surjecting onto
$(C,c)$.
\end{proof}
\begin{De}[$\mathbf{M}(R)$]
Let $R$ be a ring. The category $\mathbf{M}(R)$ is the category
of all finitely generated right $R$-modules.
\end{De}
\begin{Co}
\label{lift} Let R be a ring, let X, Y, Z and W be left flat
R-bimodules and let $\text{F}$ be one of the functors $\text{Id}$,
$\text{F}_X$, $\text{F}_{X,Y}$ and $\text{F}_{X,Y,Z,W}$. Let
$(M,m)$ be an object in $\NIL \big(\mathbf{M}(R);\text{F}\big)$.
There exists an object $(P,p)$ in the subcategory
$\NIL\big(\mathbf{P}(R);\text{F}\big)$ admitting a surjection
onto $(M,m)$.
\end{Co}
\begin{proof}
One easily checks that F has the required properties.
\end{proof}

In the rest of this section we will use the following conventions: Let $R$ be a ring and let $X$, $Y$, $Z$ and $W$ be left flat
$R$-bimodules. Let $s$ be an element of the center of $R$ which is
not a zero divisor and satisfies $s \cdot x = x \cdot s$ for all
elements $x \in X$ and similar conditions for $Y$, $Z$ and $W$.
Above we defined the functor F$\colon \Mod (R) \ra \Mod (R)$ to be one of the functors Id, F$_X$, F$_{X,Y}$ or F$_{X,Y,Z,W}$. We define F$_s \colon \Mod (R_s) \ra \Mod (R_s)$ to be respectively  Id, F$_{{}_s X_s}$, F$_{{}_s X_s,{}_s Y_s}$ or F$_{{}_s X_s,{}_s Y_s,{}_s Z_s,{}_s W_s}$. 
\begin{De}
\begin{enumerate}\label{Imd1}
\item Let $\mathbf{P}^{\Image}(R_s)$ be the full exact
subcategory of $\mathbf{P}(R_s)$ consisting of those objects
isomorphic to $P \ot_R R_s$ for some $P \in \textbf{P}(R)$.
\item Let $\mathbf{P}^{d1}(R)$ be the exact category of finitely
generated $R$-modules of projective dimension smaller or equal to
1 with the additional assumption that $P \ot_R R_s \in
\mathbf{P}^{\Image}(R_s)$ for all objects $P$ in
$\mathbf{P}^{d1}(R)$.
\end{enumerate}
\end{De}
For all full subcategories $\mathcal{C}$ of $\Mod (R)$ we define the exact functor 
\begin{align*}
\text{S} \colon \NIL(\mathcal{C};\text{F}) & \ra \NIL(\mathcal{C};\text{F}) \\
(M,m) & \mapsto(M,m \cdot s). 
\end{align*}
Let $\omega$ be the category where objects are non-negative integers and there is exactly one morphism from $i$ to $j$ if and only if $i$ is smaller or equal to $j$. The category $\omega$ is small and filtering.

In the following we will consider colimits of
functors from $\omega$ to the category of exact categories
$\mathcal{E}$ of the following type. Define
a functor from $\omega$ to $\mathcal{E}$ by sending every object of $\omega$ to a fixed category $\NIL$ and the morphisms from $i$ to
$j$ to the $(i-j)$-fold composition of S with its self. A colimit in this case consist of an exact category
$\colim \NIL$ and structure functors
$$
\text{S}_j \colon \NIL \ra \colim \NIL
$$
which are universal among functors making the diagram
$$
\xym{\cdots \ar[r] & \NIL \ar[dr]_-{\text{S}_{j-1}}
\ar[r]^{\text{S}} & \NIL \ar[d]^-{\text{S}_{j}}
\ar[r]^{\text{S}}  & \NIL
\ar[dl]^-{\text{S}_{j+1}} \ar[r]& \cdots \\
&& \colim \NIL && }
$$
commutative.  It is a result of Quillen that colimits over small and filtering categories exists in the category of exact categories \cite[page 104]{Qu}.

Localization at $s$ induces a functor from $\Mod (R)$ to $\Mod (R_s)$. By the reasoning given in the first section this induces a functor
\[
\tilde{\text{L}} \colon \NIL(\mathbf{P}^{d1}(R); \text{F}) \ra \NIL(\mathbf{P}^{\Image}(R_s);\text{F}_s).
\]
Since $\tilde{L} \circ \text{S} \cong \text{S} \circ \tilde{L}$ the universal property of a colimit implies that we obtain a functor
$$
\text{L} \colon \colim \NIL(\mathbf{P}^{d1}(R); \text{F}) \ra \colim \NIL(\mathbf{P}^{\Image}(R_s);\text{F}_s).
$$
Observe that the category $\colim \NIL(\mathbf{P}^{\Image}(R_s);\text{F}_s)$ is equivalent to the category $\NIL(\mathbf{P}^{\Image}(R_s);\text{F}_s)$ since multiplication by $s$ is invertible on $R_s$.

Let $\mathcal{C}$ be an exact category. The
category $\mathcal{C}$ can be seen as a
Waldhausen category, by defining a map to be a weak equivalence
if it is an isomorphism and to be a cofibration if it is an
admissible monomorphism. Denote by $\text{\itshape{iso}}
\mathcal{C}$ the category
of weak equivalences and by $\text{\itshape{co}} \mathcal{C}$ the category of cofibrations.
If not stated otherwise exact categories are always
considered with this Waldhausen category structure.

Let $\text{L} \colon \mathcal{C} \ra \mathcal{D}$ be an exact functor between exact categories. We give $\mathcal{C}$ a new Waldhausen structure by
defining a map to be an weak equivalence if it maps to an isomorphism under $\text{L}$. Denote this new category of weak
equivalences by $\text{\itshape{w}} \mathcal{C}$. The functor $\text{L}$ induces functors
$$
\text{L}_r \colon \text{\itshape{co}} S_r \mathcal{C} \cap \text{\itshape{w}} S_r
\mathcal{C} \ra \text{\itshape{iso}} S_r \mathcal{D}
$$
for every $r \in \nn$. We call morphisms in this category \emph{trivial cofibrations}. For an object $D \in \text{\itshape{iso}} S_r \mathcal{D}$ we denote the over-category by $\text{L}_r \downarrow D$.
\begin{Le} \label{L1}
For every $M \in \colim \NIL \big( \mathbf{P}^{\Image}(R);\text{F}_s \big)$ there exists an object $P_M$ of the subcategory $\colim \NIL(\mathbf{P}(R);\text{F})$ in $\colim \NIL(\mathbf{P}^{d1}(R);\text{F})$ and an isomorphism $f_M \colon \text{L}(P_M) \ra M$.
\end{Le}
\begin{proof}
 The lemma follows since for every finitely generated $R$-module $M$ every morphism $m \colon M \ot R_s \ra \text{F}_s(M \ot R_s) $ is induced, after multiplying with $s$ sufficiently often, by an $R$-module homomorphism from $M$ to F$(M)$. Note that here we need the assumption that $s \cdot x = x \cdot s$ for all $x \in X$ and similar assumptions for $Y$, $Z$ and $W$. 
\end{proof}
\begin{Le}\label{L2}
Let $f \colon \text{L}(C) \ra M$ be an object of $\text{L}_1 \downarrow M$. By the preceding lemma we have an isomorphism $f_M \colon \text{L}(P_M) \ra M$. There exists a trivial cofibration $h \colon P_M \ra P_M$ and a morphism $g \in mor (\text{L}_1 \downarrow M)$ from $f_M \circ \text{L}(h) \colon \text{L}(P_M) \ra M$ to $f \colon \text{L}( C ) \ra M$.
\end{Le}
\begin{proof}
 Since the category $\omega$ is filtering and $S_j \circ \tilde{L} = L \circ S_j$ we can find an object $k$ in $\omega$ such that the inverse image of the diagram
 \[
\xym{   &  \text{L}(P_M) \ar[d]^-{f_M}  \\
      \text{L}(C) \ar[r]^f & M   }
 \]
 under the $k$-th structure functor is a diagram
\[
\xym{ & \tilde{\text{L}}\big((\tilde{P}_M,\tilde{p}_m)\big) \ar[d]^-{\tilde{f}_M}  \\
      \tilde{\text{L}}\big((\tilde{C},\tilde{c})\big) \ar[r]^-{\tilde{f}} & (\tilde{M},\tilde{m})}
\]
 where $(\tilde{M},\tilde{m})$ is an object in $\NIL(\mathbf{P}^{\Image}(R_s);\text{F}_s)$ and $\xym{\tilde{\text{L}}\big((\tilde{C},\tilde{c})\big) \ar[r]^-{\tilde{f}} & (\tilde{M},\tilde{m})}$ and  $\xym{\tilde{\text{L}}\big((\tilde{P}_M,\tilde{p}_m)\big) \ar[r]^-{\tilde{f}_M} & (\tilde{M},\tilde{m})}$ are in $\tilde{\text{L}}_1 \downarrow (\tilde{M},\tilde{m})$. 

Consider the following diagram of $R$-modules
\[ 
\xym{ & & (\tilde{P}_M,\tilde{p}_m) \ar[d]^{\iota_M } \\ 
  & &\tilde{\text{L}}\big((\tilde{P}_M,\tilde{p}_m)\big) \ar[d]^-{\tilde{f}_M}  \\
      (\tilde{C},\tilde{c}) \ar[r]^-{\iota} & \tilde{\text{L}}\big((\tilde{C},\tilde{c})\big)  \ar[r]^-{\tilde{f}} & (\tilde{M},\tilde{m})}
\]
where $\iota$ and $\iota_M$ are induced by the obvious $R$-module maps from $\tilde{C}$ to $\tilde{C} \ot R_s$ and from $\tilde{P}_M$ to $\tilde{P}_M \ot R_s$. In the following we will prove that we can find an $i \in \nn$ such that the image of $\tilde{f}_M \circ \iota_M \cdot s^i$ is contained in the image of $\tilde{f} \circ \iota$. Choose a generating set $p_1, \ldots , p_l$ for $\tilde{P}_M$. We have $\tilde{f}_M \circ \iota_M (p_1) = \tilde{f}(x)$ for some $x \in \tilde{C} \ot R_s$. Thus we can find an $i' \in \nn$ such that $\tilde{f}_M \circ \iota_M (p_1) \cdot s^{i'} = \tilde{f}(x' \ot \text{id})$ for some $x' \in C$. Iterated use of this argument shows that  we can find an $i \in \nn$ such that the image of $\tilde{f}_M \circ \iota_M \cdot s^i$ is contained in the image of $\tilde{f} \circ \iota$. Since $\tilde{P}_M$ is projective we obtain an $R$-module map $\tilde{g}$ making the following diagram commute.
\[ 
\xym{ & & \tilde{P}_M \ar[d]^{\iota_M } \ar[ddll]_{\tilde{g}}\\ 
  & & \tilde{P}_M \ot R_s \ar[d]^-{\tilde{f}_M \cdot s^i}  \\
      \tilde{C} \ar[r]^-{\iota} & \tilde{C} \ot R_s \ar[r]^-{\tilde{f}} & \tilde{M}.}
\]
Note that also the diagram 
\begin{align}
\begin{split}
\xym{ 
  &\tilde{\text{L}}\big((\tilde{P}_M,\tilde{p}_m)\big) \ar[d]^-{\tilde{f}_M \cdot s^i}  \ar[dl]_{\tilde{g} \ot \text{id} } \\
  \tilde{\text{L}}\big((\tilde{C},\tilde{c})\big)  \ar[r]^-{\tilde{f}} & (\tilde{M},\tilde{m})}
\end{split}
\end{align}
commutes since $\tilde{f}$ is an isomorphisms. 

The morphism claimed in the lemma is induced by $\tilde{g} \cdot s^{j}$ for an sufficient large $j$. To prove $\tilde{g} \cdot s^{j} \in mor(\tilde{\text{L}}_1 \downarrow (\tilde{M},\tilde{m}))$ we need to prove three things.

First of all we need to show that $\tilde{g}$ defines a morphism in $S_1 \NIL(\mathbf{P}^{d1}(R);\text{F}))$, i.e. that the following diagram commutes:
\[
\xym{\tilde{P}_M \ar[r]^{\tilde{p}_m} \ar[d]_{\tilde{g}} & \text{F}(\tilde{P}_M) \ar[d]^{\text{F}(\tilde{g})} \\ \tilde{C} \ar[r]^{\tilde{c}} & \text{F}(\tilde{C}).}
\]
Since diagram 1.) commutes we have 
\[
\text{F}(\iota) \circ \tilde{c} \circ \tilde{g}(p_1) = \text{F}(\iota) \circ \text{F}(\tilde{g}) \circ \tilde{p}_m(p_1).
\]
Thus 
\[
\tilde{c} \circ \tilde{g}(p_1) \cdot s^{j'} = \text{F}(\tilde{g}) \circ \tilde{p}_m(p_1) \cdot s^{j'}
\]
for some $j' \in \nn$. Iterated use of this argument proves that we can find an $j \in \nn$ such that the diagram where $\tilde{g}$ is replaced by $\tilde{g} \cdot s^j$ commutes.

The second thing to prove is that $\tilde{g} \cdot s^j $ is a trivial cofibration. The map $\tilde{g} \cdot s^j$ is in $\text{\itshape{w}} S_1 \NIL(\mathbf{P}^{d1}(R);\text{F})$ since diagram 1.) commutes and therefore $\tilde{g} \cdot s^j$ becomes an isomorphism after localization at $s$. It is a monomorphism since $\tilde{f}_M \circ \iota_M \cdot s^{i+j}$ is a monomorphism. To prove that $\tilde{g} \cdot s^j$ is a cofibration in $S_1 \NIL(\mathbf{P}^{d1}(R);\text{F})$ it remains to show that the cokernel of $\tilde{g} \cdot s^j$ is in $\mathbf{P}^{d1}(R)$. The cokernel is finitely generated and of projective dimension smaller or equal to one since we have an exact sequence
\[
\xym{0 \ar[r] & \tilde{P}_M \ar[r]^{\tilde{g} \cdot s^j} & \tilde{C} \ar[r] & \Coker (\tilde{g} \cdot s^j) \ar[r] & 0}\]
where $\tilde{P}_M$ is finitely generated and projective and $\tilde{C}$ is finitely generated and of projective dimension smaller or equal to one. Since the module $\Coker (\tilde{g} \cdot s^j)$ is an $s$-primary torsion module we have $\Coker (\tilde{g} \cdot s^j) \ot R_s \in \mathbf{P}^{\Image}(R_s)$. Thus $\tilde{g} \cdot s^j$ is a cofibration.

Thirdly we have to verify that the diagram
\[
\xym@C=6em{&  & \text{F}_s(\tilde{P}_M \ot R_s) \ar[dd]^(0.7){\text{F}_s({\tilde{f}_M \cdot s^{i+j}})} \ar[ddl]_{\text{F}_s(\tilde{g} \cdot s^j \ot \text{id})}          \\
     &  \tilde{P}_M \ot R_s \ar[ddl]_{\tilde{g} \cdot s^j \ot \text{id}} \ar[dd]^(0.7){\tilde{f}_M \cdot s^{i+j} } \ar[ur]^-{\tilde{p}_m \ot \text{id}} & \\
         & \text{F}_s(\tilde{C} \ot R_s) \ar[r]^-{\text{F}_s(\tilde{f}')}  &  \text{F}_s(\tilde{M})                \\ 
     \tilde{C} \ot R_s \ar[ru]^-{\tilde{c} \ot \text{id}} \ar[r]^-{\tilde{f}'}  & \tilde{M} \ar[ur]_-{\tilde{m}} &}
\]
commutes. By construction the only thing which has to be proven is 
\[
(\tilde{c} \ot \text{id}) \circ (\tilde{g}\cdot s^j \ot \text{id}) = \text{F}_s(\tilde{g} \cdot s^j \ot \text{id}) \circ (\tilde{p}_m \ot \text{id}). 
\]
But this follows by the argument given above. This implies that $\tilde{g} \cdot s^j$ is a morphism in $\tilde{\text{L}}_1 \downarrow (\tilde{M},\tilde{m})$.

Multiplication by $s^{i+j}$ becomes an isomorphism after localization and is a cofibration since $s$ is a non zero divisor. This implies that multiplication by $s^{i+j}$ is a trivial cofibration.

 The image of $s^{i+j}$ and $g \cdot s^j$ under the $k+j$-th structural functor give the desired morphisms.
\end{proof}
\begin{Le} \label{L3}
 For every $P \in S_r \colim \NIL(\mathbf{P}(R);\text{F})$ and pair of morphisms $f_1$ and $f_2 \colon P \ra C$ in $\text{L}_r \downarrow M$ the assumption L$_r(f_1) = \text{L}_r(f_2)$ implies that there is trivial cofibration $h$ such that $f_1 \circ h = f_2 \circ h$.
\end{Le}
\begin{proof}
The morphism $h$ is induced by multiplication by $s^i$ for sufficiently large $i$. The existence of such an $i$ follows since projective $R$-modules are $s$-torsion free.
\end{proof}
\begin{Le} \label{L5}
Assume we have a commutative diagram
$$
\xym{0 \ar[r] &\text{L}(C_1) \ar@{>->}[r] \ar[d]^{\cong} &\text{L}(C_2) \ar@{->>}[r] \ar[d]^{\cong} & \text{L}(C_3) \ar[d]^{\cong} \ar[r] & 0 \\
0 \ar[r] & D_1 \ar@{>->}[r] & D_2 \ar@{->>}[r] & D_3 \ar[r] & 0 }
$$
with exact rows and the upper row is induced by an exact sequence in the category $\colim \NIL(\mathbf{P}^{d1}(R);\text{F})$. By Lemmas \ref{7} and \ref{L2} we obtain a commutative diagram
$$
\xym{ & & \text{L}(P_1) \ar@{>->}[r] \ar[dd] \ar[dll] & \text{L}(P_2) \ar@{->>}[r] \ar[dd] & \text{F}(P_3) \ar[dd] \ar[dll]  \\
\text{L}(C_1) \ar[drr] \ar@{>->}[r] & \text{L}(C_2) \ar[drr] \ar@{->>}[r] & \text{L}(C_3) \ar[drr] & &         \\
     & & D_1  \ar@{>->}[r] & D_2 \ar@{->>}[r] & D_3.}
$$
with $P_1$, $P_2$ and $P_3 \in \NIL(\mathbf{P}(R);\text{F}
)$. We can find a map from $P_2$ to $C_2$ making the given diagram commutative. 
\end{Le}
\begin{proof}
The result follows in a similar way as Lemma \ref{L2} is proven.
\end{proof}
\subsection{A Long Exact Sequence}
In this section we develop a long exact sequence relating the $K$-groups of certain exact categories.
\begin{Pro} \label{locle} Let $\mathcal{C}$ and $\mathcal{D}$ be exact categories and let $\text{L} \colon \mathcal{C} \ra \mathcal{D}$ be an exact
functor. We give $\mathcal{C}$ a new Waldhausen structure by
defining a map to be a weak equivalence if it maps to an isomorphism under $\text{L}$. Denote this new category of weak
equivalences by $\text{\itshape{w}} \mathcal{C}$. Let
$\Ker(\text{L})$ be the kernel of $\text{L}$, i.e., the full
Waldhausen subcategory of all objects $C$ in $\mathcal{C}$ such
that there exists an isomorphism $\text{L}(C) \cong 0$ in $\mathcal{D}$.

If the functor $\text{L}$ induces a homotopy equivalence
$$
\text{\itshape{co}} S_r \mathcal{C} \cap \text{\itshape{w}} S_r
\mathcal{C} \ra \text{\itshape{iso}} S_r \mathcal{D}
$$
after taking the nerve and realization for every $r \in \nn$, then
$$
\xym{\Ker(\text{L}) \ar[r] & \mathcal{C} \ar[r]^{\text{G}} &
\mathcal{D}}
$$
induces a long exact sequence on $K$-groups.
\end{Pro}
\begin{proof} Consider the inclusion functor
$$
\text{I} \colon \Ker (\text{L}) \ra \mathcal{C}.
$$
Corollary I.1.5.7. and Corollary (2) of Lemma I.1.4.1 in \cite{Wa2} gives that
 $$
 \xym{ \big|s_{\bullet} \Ker (\text{L}) \big| \ar[r]^-{\text{I}} &
 \big| s_{\bullet} \mathcal{C} \big| \ar[r] &
 \big| s_{\bullet} S_{\bullet} \big(\Ker (\text{G}) \ar[r]^-{\text{I}} &
 \mathcal{C}\big) \big| }
 $$
 is a fibration up to homotopy. The main step is to prove that $\text{L}$
 identifies
 $$\xym{\big| s_{\bullet} S_{\bullet} \big(\Ker (\text{L}) \ar[r]^-{\text{I}} &
 \mathcal{C}\big) \big|}$$
 with
 $$
 \big| N_{\bullet} \text{\itshape{iso}} S_{\bullet} \mathcal{D}\big|,
 $$
 where $N_{\bullet}$ is the nerve of a category.
 
Since $\Ker (\text{L})$ is a Waldhausen subcategory of  $\mathcal{C}$, the category $S_m \Ker (\text{L})$  is a Waldhausen subcategory of $S_m \mathcal{C}$.
 Thus following Waldhausen \cite[page 344]{Wa2}, we can replace
\[
\xym{s_{\bullet} S_\bullet (\Ker (\text{L}) \ar[r]^-{\text{I}} & \mathcal{C} \big)}
\]
  by
\[
 \xym{s_\bullet F_\bullet \big( \mathcal{C}, \Ker (\text{L})\big).}
\]
Observe that by reversal of priorities we can replace
 the bisimplicial set
 $$
\xym{(m,n) \ar@{|->}[r] & s_m F_n\big(\mathcal{C}, \Ker (\text{L})\big)}
$$
 by the equivalent bisimplicial set
 $$
 \xym{(m,n) \ar@{|->}[r] & obj\big(F_n \big( S_m \mathcal{C}, S_m \Ker (\text{L}) \big)\big).}
 $$
 The functor $\text{L}$ induces a map of bisimplicial sets from
 $$
 \xym{(m,n) \ar@{|->}[r] & obj\big(F_n \big(S_m \mathcal{C}, S_m \Ker(\text{L})
 \big)\big)}
 $$
 to
 $$
 \xym{(m,n) \ar@{|->}[r] & N_n \big( \text{\itshape{iso}} S_m \mathcal{D} \big).}
 $$

 The next step is to identify $obj \big( F_n \big(S_m \mathcal{C}, S_m \Ker (\text{L}) \big)\big)$ also with
 the nerve of a category. The bisimplicial set
 $$
 \xym{(m,n) \ar@{|->}[r] & obj \big( F_n \big(S_m \mathcal{C}, S_m \Ker (\text{L}) \big)\big)}
 $$
 is equivalent to
 $$
 \xym{(m,n) \ar@{|->}[r] & N_n ( \text{\itshape{co}} S_m \mathcal{C}
 \cap \text{\itshape{w}} S_m \mathcal{C}).}
 $$
 Thus by the realization lemma, to prove the statement it suffices
 to show that for each $r \ge 0$ the induced functor
 \begin{center}
 \makebox[0pt]{
 $$\xym{\text{L}_r \colon \text{\itshape{co}} S_r \mathcal{C}
 \cap \text{\itshape{w}} S_r \mathcal{C}
 \ar[r] & \text{\itshape{iso}} S_r \mathcal{D}}$$}
 \end{center}
 realizes to a homotopy equivalence. But this is one of our
 assumptions.
\end{proof}
\begin{Pro}\label{Pro.2}
Let the notation be as in the preceding proposition. The functor 
\begin{center}
 \makebox[0pt]{
 $$\xym{\text{L}_r \colon \text{\itshape{co}} S_r \mathcal{C}
 \cap \text{\itshape{w}} S_r \mathcal{C}
 \ar[r] & \text{\itshape{iso}} S_r \mathcal{D}}$$}
\end{center}
induces a homotopy equivalent after taking the nerve and realization for every $r \in \nn$ if there exists a class $\mathcal{P}$ of objects of $\mathcal{C}$ with the following properties:
\begin{enumerate}
\item For every object $D$ of the category $\mathcal{D}$ there exists $P_D \in \mathcal{P}$ and an isomorphism $f_D \colon \text{L}(P_D) \ra D$ in $\mathcal{D}$.
\item For every $(f \colon \text{L}(C) \ra D) \in obj(\text{L}_1 \downarrow D)$ there exists a trivial cofibration $h$ from $P_D$ to $P_D$ and a morphism $g \in \text{L}_1 \downarrow D$ from $f_D \circ \text{L}(h) \colon \text{L}(P_D) \ra D$ to $f \colon \text{L}(C) \ra D$.
\item For every $P \in S_r \mathcal{P}$ and pair of morphisms $f_1, f_2 \colon P \ra C$ in $\text{L}_r \downarrow D$ the assumption L${}_r(f_1) = \text{L}_r(f_2)$ implies that there is a morphism $h \colon C' \ra P$ in $\text{L}_r \downarrow D$ such that $f_1 \circ h = f_2 \circ h$.
\item Let 
$$
\xym{0 \ar[r] & D_1 \ar@{>->}[r] & D_2 \ar@{->>}[r] & D_3 \ar[r] & 0}
$$
be an exact sequence in $\mathcal{D}$. Let $P_1, P_3 \in \mathcal{P}$ such that $\text{L}(P_1) \cong D_1$ and $\text{L}(P_3) \cong D_3.$
There exists $P_2 \in \mathcal{P}$ with $\text{L}(P_2) \cong D_2$ and there is a commutative diagram 
$$
\xym{0 \ar[r] & \text{L}(P_1) \ar[d]^-{\cong} \ar@{>->}[r] & \text{L}(P_2) \ar[d]^-{\cong} \ar@{->>}[r] & \text{L}(P_3) \ar[d]^-{\cong} \ar[r] & 0 \\
0 \ar[r] & D_1 \ar@{>->}[r] & D_2 \ar@{->>}[r] & D_3 \ar[r] & 0}
$$
where the upper row is induced by an exact sequence in $\mathcal{C}$.
\item 
Assume we have a commutative diagram
$$
\xym{0 \ar[r] &\text{L}(C_1) \ar@{>->}[r] \ar[d]^{\cong} &\text{L}(C_2) \ar@{->>}[r] \ar[d]^{\cong} & \text{L}(C_3) \ar[d]^{\cong} \ar[r] & 0 \\
0 \ar[r] & D_1 \ar@{>->}[r] & D_2 \ar@{->>}[r] & D_3 \ar[r] & 0 }
$$
with exact rows and the upper row is induced by an exact sequence in $\mathcal{C}$. By property 2. and 4. we obtain a commutative diagram
$$
\xym{ & & \text{L}(P_1) \ar@{>->}[r] \ar[dd] \ar[dll] & \text{L}(P_2) \ar@{->>}[r] \ar[dd] & \text{L}(P_3) \ar[dd] \ar[dll]  \\
\text{L}(C_1) \ar[drr] \ar@{>->}[r] & \text{L}(C_2) \ar[drr] \ar@{->>}[r] & \text{L}(C_3) \ar[drr] & &         \\
     & & D_1  \ar@{>->}[r] & D_2 \ar@{->>}[r] & D_3.}
$$
with $P_1$, $P_2$ and $P_3 \in \mathcal{P}$. We can find a trivial cofibration from $P_2$ to $C_2$ making the given diagram commutative. 
\end{enumerate}
\end{Pro}
\begin{proof}
The result will follow from Quillens Theorem A once we show that $\text{G}_r \downarrow D$ is nonempty and cofiltering and therefore contractible for every $D \in \text{\itshape{iso}} S_r \mathcal{D}$.

For $r=0$ there is nothing to prove. For $r=1$ Condition 1. implies that $\text{L}_1 \downarrow D$ is nonempty.

Given two objects 
$$
\xym{\text{L}(C) \ar[r]^-{f} & D & \ar[l]_-{f'} \text{L}(C')}
$$
in the category $\text{L}_1 \downarrow D$. By Condition 2. we have an object $P$ and morphisms $h$ and $h'$ such that the following diagram commutes.
$$
\xym{                        & \text{L}(P) \ar[d]^-{\text{L}(h')} \ar[dddl]  &   \\ 
                             & \text{L}(P) \ar[d]^-{\text{L}(h)}  \ar[ddr] &   \\
                             & \text{L}(P) \ar[d]        &   \\
     \text{L}(C) \ar[r]^-{f} &  D                  &  \text{L}(C'). \ar[l]_-{f'}}
$$
Thus there exist an object in L$_1 \downarrow D$ which maps on both objects.

Given two morphisms 
\begin{center}
 \makebox[0pt]{
 $\xym@C=1.5em{\text{L}(f_1) \colon (\text{L}(C) \ar[r]^-{f} &
 D) \ar[r] & (\text{L}(C') \ar[r]^-{f'} & D)}
 $}
\end{center}
 \begin{center}
 \makebox[0pt]{
 $\xym@C=1.5em{\text{L}(f_2) \colon (\text{L}(C)
 \ar[r]^-{f} & D )
 \ar[r] & (\text{L}(C') \ar[r]^-{f'} & D)}$}
 \end{center}
in the category $\text{L}_1 \downarrow D$. The morphisms $f' \circ \text{L}(f_1)$ and $f' \circ \text{L}(f_2)$ define objects in L$_1 \downarrow D$. By the reasoning given above and since maps in $\text{\itshape{iso}} \mathcal{D}$ are isomorphisms there is an object $P_D \in \mathcal{P}$ and a map $g \colon P_D \ra C$ such that $\text{L}(f_1 \circ g)  = \text{L}(f_2 \circ g)$. By Condition 3. we have a morphism such that $f_1 \circ g \circ h = f_2 \circ g \circ h$. Thus $\text{L}_1 \downarrow D$ is non empty and cofiltering.

For $r \ge 2$, the arguments can be extended in the following
 manner. An object of $\text{L}_r \downarrow D$ amounts to a diagram
 $$
C =\xym{0 \, \, \ar@{>->}[r] & C^{1,0}
 \ar@{->>}[d] \, \,\ar@{>->}[r] & \ldots \, \, \ar@{>->}[r]&
 C^{r,0}
 \ar@{->>}[d] \\ & 0 \, \, \ar@{>->}[r] & \ldots \, \, \ar@{>->}[r] &
 C^{r,1} \ar@{->>}[d] \\ & & & \ldots
 \ar@{->>}[d] \\ & & 0 \, \, \ar@{>->}[r] & C^{r,r-1}
 \ar@{->>}[d] \\ & & & 0}$$
in $S_r \mathcal{C}$ and a map $f \in \text{\itshape{iso}} S_r \mathcal{D}$ from L$_r(C)$ to the diagram
 $$
 D = \xym{0 \, \, \ar@{>->}[r] & D^{1,0} \ar@{->>}[d]
 \, \, \ar@{>->}[r] & \ldots \, \, \ar@{>->}[r]&
 D^{r,0}
 \ar@{->>}[d] \\ & 0 \, \, \ar@{>->}[r] & \ldots \, \, \ar@{>->}[r] &
 D^{r,1} \ar@{->>}[d] \\ & & & \ldots
 \ar@{->>}[d] \\ & & 0 \, \, \ar@{>->}[r]& D^{r,r-1}
 \ar@{->>}[d] \\ & & & 0.}
 $$
By Condition 1. we can find $P^{i,i-1} \in \mathcal{P}$ such that $\text{G}(P^{i,i-1}) \cong D^{i,i-1}$ and by Condition 4. we can find objects $P^{i+1,i-1} \in \mathcal{P}$ and morphisms such that we obtain a commutative diagram 
$$
\xym{\text{L}(P^{i,i-1}) \ar[d]^{\cong} \ar@{>->}[r] & \text{L}(P^{i+1,i-1}) \ar[d]^{\cong} \ar@{->>}[r] & \text{L}(P^{i+1,i}) \ar[d]^{\cong} \\
D^{i,i-1} \ar@{>->}[r] & D^{i+1,i-1} \ar@{->>}[r] & D^{i+1,i} }
$$
where the rows are exact.
Applying this construction $r$-times we get an object in $\text{L}_r \downarrow D$. Thus $\text{L}_r \downarrow D$ is non empty. 

Suppose we have two objects
$$
\xym{\text{L}_r(C) \ar[r] & D & \text{L}_r(C') \ar[l]}
$$
in the category $\text{L}_r \downarrow D$. The components of $C$ and $C'$ are denoted $C^{i,j}$ and $C'^{i,j}$ respectively. As in the proof when $r=1$ we can construct for every part
$$
\xym{\text{L}(C^{i,i-1}) \ar[r] & D^{i,i-1} & \text{L}(C'^{i,i-1}) \ar[l]}
$$
of the diagram an object $P^{i,i-1} \in \mathcal{P}$ and morphisms such that 
$$
\xym{ & \text{L}(P^{i,i-1}) \ar[ld] \ar[rd] \ar[d] & \\
     \text{L}(C^{i,i-1}) \ar[r] & D^{i,i-1} & \text{L}(C'^{i,i-1}) \ar[l]}
$$
commutes. Using Condition 4. we obtain a commutative diagram
$$
\xym{\text{L}(P^{i,i-1}) \ar[d]^{\cong} \ar@{>->}[r] & \text{L}(P^{i+1,i-1}) \ar[d]^{\cong} \ar@{->>}[r] & \text{L}(P^{i+1,i}) \ar[d]^{\cong} \\
D^{i,i-1} \ar@{>->}[r] & D^{i+1,i-1} \ar@{->>}[r] & D^{i+1,i} }
$$
where the rows are exact. Condition 5. gives morphisms from $\text{L}(P^{i+1,i-1})$ to $\text{L}(C^{i+1,i-1})$. Applying this construction $r$-times we get an element which maps onto the two objects. 

Given two morphisms
\begin{center}
 \makebox[0pt]{
 $\xym@C=1.5em{\text{L}_r(f_1) \colon (\text{L}_r(C) \ar[r]^-{f} &
 D) \ar[r] & (\text{L}_r(C') \ar[r]^-{f'} & D)}
 $}
\end{center}
 \begin{center}
 \makebox[0pt]{
 $\xym@C=1.5em{\text{L}_r(f_2) \colon (\text{L}_r(C)
 \ar[r]^-{f} & D )
 \ar[r] & (\text{L}_r(C') \ar[r]^-{f'} & D)}$}
 \end{center}
in the category $\text{L}_r \downarrow D$. As above the morphisms $f' \circ \text{L}_r(f_1)$ and $f' \circ \text{L}_r(f_2)$ define objects in L$_r \downarrow D$. By the reasoning given above and since maps in $\text{\itshape{iso}} S_r \mathcal{D}$ are isomorphisms there is an object $P_D \in S_r \mathcal{P}$ and a map $g \colon P_D \ra C$ such that $\text{L}_r(f_1 \circ g)  = \text{L}_r(f_2 \circ g)$. By Condition 3. we have a morphism such that $f_1 \circ g \circ h = f_2 \circ g \circ h$. Thus $\text{L}_r \downarrow D$ is non empty and cofiltering.
\end{proof}
\subsection{The Long Exact Localization Sequence}
We apply the results of the preceding sections to obtain a long exact localization sequence for the $K$-groups of Nil-categories. As a corollary we obtain localization results for Nil-groups.

We define $\mathbf{H}_s(R)$ to be the exact category of finitely generated
$R$-modules with the property that modules have projective
dimension smaller or equal to~1 and are $s$-primary torsion. We obtain an inclusion-functor
\[
\text{I} \colon \colim \NIL(\mathbf{H}_s(R);\text{F}) \ra \colim \NIL(\mathbf{P}^{d1}(R);\text{F}).
\]
\begin{Th}
\label{loca} The sequence
\[
\xym{\colim \NIL \big(\mathbf{H}_s(R);\text{F} \big) \ar[r]^-{I} & \colim
\NIL\big(\mathbf{P}^{d1}(R);\text{F} \big) \ar[r]^-{L} & \colim
\NIL\big(\mathbf{P}^{\Image}(R_s); \text{F}_s \big)}
\]
of functors induces a long exact sequence on $K$-theory.
\end{Th}
\begin{proof}
The theorem follows from Proposition \ref{locle} and \ref{Pro.2}, Lemma \ref{7}, \ref{L1}, \ref{L2}, \ref{L3} and \ref{L5} and the fact that $\colim \NIL \big( \mathbf{H}_s(R);\text{F} \big)$ is the kernel of L. 
\end{proof}

In the last part of this section we derive Theorem~\ref{Th2} out
of Theorem~\ref{loca}. Note that we are again sloppy with the notation since for F$= \text{F}_{X,Y}$ or F$_{X,Y,Z,W}$ the first identity should for example be 
\[
K_i \big(\colim \NIL(\mathbf{H}_s(R);\text{F})\big) \cong
K_i\big(\mathbf{H}_s(R) \times \mathbf{H}_s(R)\big).
\]
\begin{Le}
\label{cong2} We have
$$
K_i\big(\colim \NIL(\mathbf{H}_s(R);\text{F})\big) \cong
K_i\big(\mathbf{H}_s(R)\big)
$$
for $i \geq 0$.
\end{Le}
\begin{proof}
 We have an equivalence of categories
 $$
 \colim \NIL \big(\mathbf{H}_s(R);\text{F} \big) \cong  \mathbf{H}_s(R) 
 $$
 since every $R$-module morphism with a source which is finitely generated and $s$-torsion becomes the trivial morphism after multiplying with $s$ sufficiently often.
\end{proof}
\begin{Le}
\label{cong3} We have
$$
 K_i \big(\colim \NIL (\mathbf{P}^{\Image}(R_s);\text{F}_s )
\big) \cong K_i\big( \NIL(\mathbf{P}(R_s); \text{F}_s) \big)
$$
for $i \ge 1$ and there is an injective map on $K_0$.
\end{Le}
\begin{proof}
 The categories $\colim \NIL (\mathbf{P}^{\Image}(R_s);\text{F}_s)$ and $\NIL (\mathbf{P}^{\Image}(R_s);\text{F}_s)$ are equivalent since multiplication by $s$ is invertible on $R_s$. This implies that we have
 $$
 K_i \big(  \colim \NIL(\mathbf{P}^{\Image}(R_s); \text{F}_s) \big) \cong
 K_i\big( \NIL (\mathbf{P}^{\Image}(R_s);\text{F}_s )\big).
 $$
 Since $\NIL \big(\mathbf{P}^{\Image}(R_s); \text{F}_s \big)$
 contains $\NIL(\mathbf{F}(R_s);\text{F}_s)$, where $\mathbf{F}(R_s)$ is the category of all finitely generated free $R_s$-modules,
 $\NIL \big( \mathbf{P}^{\Image}(R_s);\text{F}_s \big)$ is
 cofinal in $\NIL(\mathbf{P}(R_s);\text{F}_s)$. Therefore
 $$K_i\big( \NIL (\mathbf{P}^{\Image}(R_s);\text{F}_s) \big) \cong
 K_i \big( \NIL(\mathbf{P}(R_s);\text{F}_s) \big)$$
 for $i \ge 1$ and there is an injective map on $K_0$.
\end{proof}
\begin{Le}
\label{2} Let $(M,m)$ be an object in $\NIL
\big(\mathbf{P}^{d1}(R);\text{F} \big)$ and let $(P,p)$ be an
object in the subcategory $\NIL(\mathbf{P}(R);\text{F})$. Let
$$
f \colon (P,p) \ra (M,m)
$$
be a surjective morphism in $\NIL \big(\mathbf{P}^{d1}(R); \text{F}\big)$. The tuple $(\Ker(f),p|_{\Ker(f)})$ is a well defined object in $\NIL(\mathbf{P}(R);\text{F})$.
\end{Le}
\begin{proof}
 The exact categories $\NIL \big(\mathbf{P}^{d1}(R);\text{F} \big)$ and $\NIL(\mathbf{P}(R); \text{F})$ are subcategories of the abelian
category $\NIL\big(\Mod (R);\text{F} \big)$. Thus
$(\Ker(f),p|_{\Ker(f)})$ is a well
defined object in $\NIL\big(\Mod (R);\text{F} \big)$. An application of Schanuel's Lemma gives that $(\Ker(f),p|_{\Ker(f)})$ is an object
 in the subcategory $\NIL(\Mod (R);\text{F})$.
\end{proof}
\begin{Le}
\label{cong}
 We have
$$
K_i\big( \NIL ( \mathbf{P}^{d1}(R);\text{F} ) \big) \cong
K_i\big( \NIL(\mathbf{P}(R);\text{F}) \big)
$$
 for $i \geq 0$.
\end{Le}
\begin{proof}
 The preceding lemma and the lift constructed in Section~\ref{Lifts}
 (Corollary~\ref{lift}) give that any object in
 $\NIL \big(\mathbf{P}^{d1}(R);\text{F}\big)$ has a length one
 $\NIL(\mathbf{P}(R);\text{F})$-resolution. To apply the
 Resolution Theorem, it is necessary to check that $\NIL(\mathbf{P}(R);\text{F})$ is closed
 under kernels in $\NIL\big(\mathbf{P}^{d1}(R);\text{F} \big)$, i.e., if
 $$
 \xym{0 \ar[r] & (M,m) \ar[r] & (P',p') \ar[r] &
 (P,p) \ar[r] & 0}
 $$
 is an exact sequence in $\NIL \big(\mathbf{P}^{d1}(R);\text{F} \big)$
 with $(P',p')$ and
 $(P,p)$ in $\NIL(\mathbf{P}(R); \text{F})$, then $(M,m)$ is also in
 $\NIL(\mathbf{P}(R);\text{F})$. This follows since the category
 $\mathbf{P}(R)$ has this property in $\mathbf{P}^{d1}(R)$.
 \end{proof}
\begin{Th}
\label{Th2}
Localization at $s$ induces a long exact sequence
\begin{align*}
\xym{ \cdots \ar[r] & K_i \big( \mathbf{H}_s(R) \big) 
\ar[r]^-I & \colim K_i \big( \NIL(\mathbf{P}(R);\text{F}) \big)
\ar[r]^-{\text{L}} & } \\
\xym{\ar[r]^-{\text{L}} & K_i \big( \NIL(\mathbf{P}(R_s); \text{F}_s)
\big) \ar[r] & K_{i-1} \big( \mathbf{H}_s(R) \big) \ar[r] & \cdots,}
\end{align*}
where $K_0 \big( \NIL(\mathbf{P}(R);\text{F}) \big) \ra K_0 \big(
\NIL(\mathbf{P}(R_s);\text{F}_s) \big)$ is not necessarily surjective.
\end{Th}
\begin{proof}
Combining the Lemmas~\ref{cong2}, \ref{cong3} and \ref{cong}, we
get that Theorem~\ref{Th2} is implied by Theorem~\ref{loca}.
\end{proof}
\begin{Co}
\label{Th} Let $R$ be a ring and let X, Y, Z and W be left flat
$R$-bimodules and let $s$ be a central non zero divisor which satisfies $s \cdot x = x \cdot s$ for all
$x \in X$ and similar conditions for $Y$, $Z$ and $W$. We define the functor F$\colon \Mod (R) \ra \Mod (R)$ to be one of the functors Id, F$_X$, F$_{X,Y}$ or F$_{X,Y,Z,W}$ and F$_s \colon \Mod (R_s) \ra \Mod (R_s)$ to be respectively  Id, F$_{{}_s X_s}$, F$_{{}_s X_s,{}_s Y_s}$ or F$_{{}_s X_s,{}_s Y_s,{}_s Z_s,{}_s W_s}$.We obtain an isomorphism
$$ \zz [t,t^{-1}] \ot_{\zz[t]} \Nil_i(\mathbf{P}(R);\text{F}) \cong
\Nil_i(\mathbf{P}(R_s); \text{F}_s),$$
for all $i \in \zz$, and $t$ acts on $\Nil_i(\mathbf{P}(R);\text{F})$ via the
map induced by the functor S.
\end{Co}
\begin{proof}
From now on, $\zz[t,t^{-1}] \ot_{\zz[t]} \Nil_i(\mathbf{P}(R);\text{F})$ is denoted $\Nil_i(\mathbf{P}(R);\text{F})_s$. Basic algebra shows that $\Nil_i(\mathbf{P}(R);\text{F})_s$ is isomorphic to $\colim \Nil_i(\mathbf{P}(R);\text{F})$.
We have $K_i\big(\NIL(\mathbf{P}(R);\text{F})\big) = \Nil_i(\mathbf{P}(R);\text{F}) \op K_i(R)$. The map which is induced by the functor S respects the given direct sum decomposition and is the identity on $K_i(R)$. Thus 
$$
\Nil_i(R;\text{F})_s = \Ker \Big( \colim K_i
\big(\NIL(\mathbf{P}(R);\text{F})\big) \ra K_i(R) \Big).
$$
In Theorem~\ref{Th2} it is proven that the following commutative
diagram has an exact row in the middle. The lowest row is the
localization sequence of algebraic $K$-theory and therefore exact.
\begin{center}
\makebox[0pt]{%
$$\xym@=1em{\ar[r]  & 0 \ar[r] \ar[d] &
\Nil_i(\mathbf{P}(R);\text{F})_s \ar[r] \ar[d] & \Nil_i(\mathbf{P}(R_s);\text{F}_s)
\ar[r] \ar[d] & \\  \ar[r] & 
K_i\big(H_s(R)\big) \ar[r] \ar[d]^{\cong} & \colim
K_i\big(\NIL(\mathbf{P}(R);\text{F})\big) \ar[r] \ar[d] &
K_i\big((\NIL(\mathbf{P}(R_s);\text{F}_s)\big) \ar[r] \ar[d] &  \\
\ar[r] & K_i\big(\mathbf{H}_s(R)\big) \ar[r] &
K_i(R) \ar[r] & K_i(R_s) \ar[r] &  .}$$ }
\end{center}
The snake lemma implies now that  $\Nil_i(\mathbf{P}(R);\text{F})_s$ and
$\Nil_i(\mathbf{P}(R_s);\text{F}_s)$ are isomorphic, which is the
statement of Corollary~\ref{Th} for $i \geq 1$.

To obtain the isomorphism for $i = 0$ note that
$$
\Ker \Big( K_0 \big( \NIL (
\mathbf{P}^{\Image}(R_s);\text{F}_s) \big) \ra K_0 \big(
\mathbf{P}^{\Image}(R_s)\big)
\Big)
$$
and $\Nil_0(R_s;\text{F}_s)$ are isomorphic since elements of
the form $[ (P,Q,0,0) ]$ are trivial in
$\Nil_0(R_s;\text{F}_s)$.

To obtain the statement for $i < 0$ note that the suspension
construction commutes with localization, i.e. $(\Sigma R)_s =
\Sigma R_s$.
\end{proof}
\begin{Rem}
We have not used the fact that the morphisms are nilpotent. Thus Corollary~\ref{Th} stays valid if we replace $\Nil$ by $\End$.
\end{Rem}
\section{Nil-Groups as Modules over the Ring of Witt Vectors}
We develop a Witt vector-module structure on a
certain class of Nil-groups including the important cases
$\Nil_i(R G; \al)$ and $\Nil_i(RG; RG_{\al}, RG_{\beta})$ where
$G$ is a group, $R$ is a commutative ring and $\al$ and $\beta$
are inner group automorphisms.
\subsection{Nil-Groups as Modules over $\End_0$}
\label{pairing}
We define a $\End_0$-module structure on certain
Nil-groups. As an application we obtain an $\End_0(R)$-module
structure on $\Nil_i(\mathbf{P}(\La);\text{F})$ if $\La$ is an algebra over a
commutative ring $R$. On $\Nil_i(\La)$ a similar module structure
is defined by Weibel \cite{We1}. The main ingredients are exact pairings. For a definition of an exact pairing see \cite{Wa1,Wa4}. 
\begin{De}
\begin{enumerate} 
\item Let $\mathcal{B}$ be an exact category with an exact pairing
$\mathcal{B} \times \mathcal{B} \rightarrow \mathcal{B}$. We define an exact pairing
\begin{align*}
\END(\mathcal{B}) \times \END(\mathcal{B}) & \ra
\END(\mathcal{B}) \\ 
\big( (B_1,b_1), (B_2,b_2)\big) & \mapsto (B_1
\times B_2, b_1 \times b_2 ).
\end{align*}
\item 
Let $\mathcal{A}$ be an abelian category, let $\mathcal{C}
\subseteq \mathcal{A}$ be a full subcategory which is closed under
extension and let $\text{F} \colon
\mathcal{A} \ra \mathcal{A}$ be an exact functor. Assume we have an exact
pairing $\mathcal{B} \times \mathcal{A} \ra \mathcal{A}$
which restricts to $\mathcal{C}$ together with a natural transformation $U_{\mathcal{C}}$ between the functors
\begin{align*}
\mathcal{B} \times \mathcal{C} & \ra \mathcal{C} \\
\text{F}_1 \colon (B,C) & \mapsto B \times \text{F}(C)\\
\text{F}_2 \colon (B,C) & \mapsto \text{F}(B \times C).
\end{align*}
We define an exact pairing
\begin{align*}
\END(\mathcal{B}) \times \NIL(\mathcal{C};\text{F}) & \ra
\NIL(\mathcal{C};\text{F}) \\ \big( (B,b), (C,c)\big) & \mapsto (B
\times C, U_\mathcal{C}(B,C) \circ (b \times c)).
\end{align*}
\end{enumerate}
\end{De}
The proof that the given pairings are well-defined is left to the reader.
\begin{Pro}
Let the notation be as in the preceding definition. If
\[ 
\xym{\END(\mathcal{B}) \times \END(\mathcal{B}) \times 
\END(\mathcal{B}) \ar[r]^-{\big(\text{Id},(-,-)\big)} \ar[d]_-{\big((-,-),\text{Id} \big)} & \END(\mathcal{B}) \times \END(\mathcal{B}) \ar[d]^{(-,-)}\\
\END(\mathcal{B}) \times \END(\mathcal{B}) \ar[r]^{(-,-)} & \END(\mathcal{B})} 
\]
and
\[
\xym{\END(\mathcal{B}) \times \END(\mathcal{B}) \times \NIL(\mathcal{C};\text{F})  \ar[r]^-{\big(\text{Id},(-,-)\big)} \ar[d]_-{\big((-,-),\text{Id} \big)} &  \END(\mathcal{B}) \times \NIL(\mathcal{C};\text{F}) \ar[d]^-{(-,-)} \\
\END(\mathcal{B}) \times \NIL(\mathcal{C};\text{F}) \ar[r]^{(-,-)}& \NIL(\mathcal{C};\text{F}) }
\]
commutes up to natural isomorphism then
$\End_0(\mathcal{B})$ carries a ring structure and
$\Nil_i(\mathcal{C};\text{F})$ is an $\End_0(\mathcal{B})$-module.
\end{Pro}
\begin{proof}
The machinery developed by Waldhausen \cite{Wa1,Wa4} implies that
we get a $K_0\big(\END(\mathcal{B})\big)$-module structure on
$K_i\big(\NIL(\mathcal{C};\text{F})\big)$. Pairing with objects of
the form $(B,0)$ reflects $\END(\mathcal{B})$ into $\mathcal{B}$
and $\NIL(\mathcal{C};\text{F})$ into $\mathcal{C}$. Thus the
$K_0(\END(\mathcal{B}))$-module structure restricts to an
$\End_0(\mathcal{B})$-module structure on
$\Nil_i(\mathcal{C};\text{F})$.
\end{proof}
\begin{Co}\label{nilmodulestructure}
Let $\La$ be an algebra over a commutative ring $R$. Let $X$, $Y$,
$Z$ and $W$ be arbitrary left flat $\La$-bimodules. The groups
$\Nil_i(\La)$, $\Nil_i(\La;X)$, $\Nil_i(\La;X,Y)$ and
$\Nil_i(\La;X,Y,Z,W)$ are modules over the ring $\End_0(R)$ for
all $i \in \zz$.
\end{Co}
\begin{proof}
The pairings which are induced by the tensor product and the obvious natural transformations satisfy the assumptions of the preceding proposition.
\end{proof}
In the following, this module multiplication is denoted by $\ast$.
\subsection{Operations on Nil-Groups}
\label{operations}
For the whole section we assume that $\mathcal{A}$ is an abelian category, $\text{F} \colon
\mathcal{A} \ra \mathcal{A}$ is an exact functor and $\mathcal{C}
\subseteq \mathcal{A}$ is a full subcategory which is closed under
extension. Furthermore $R$ is a ring, $G$ is a group, $X$ and $Y$ are arbitrary $RG$-bimodules and $\al$ and $\beta$ are inner group automorphisms induced by group elements $g$ and $g'$. Recall the definition of the $RG$-bimodule $RG_{\al}$ given in Remark \ref{Rem}.
\subsubsection{Frobenius Operations on Nil-Groups}
On Bass Nil-groups, for a natural number $n$ the $n$-th Frobenius
is defined to be the map induced by the functor whose value at
$(P,p)$ is $(P,p^n)$. The main problem with the definition of
the Frobenius operation on more general Nil-groups is that in
general $\text{F}(c) \circ c$ is an object in
$\NIL(\mathcal{C}; \text{F}^2)$ and not in
$\NIL(\mathcal{C};\text{F})$. To get around this problem, we use
the assumption that we have a natural transformation between
$\text{F}^\ell$ and the identity. 
\begin{De}[Frobenius]
\begin{enumerate}
\item let $\text{F}^0$
be the identity functor $\mathcal{C} \ra \mathcal{C}$ and for a
natural number $n \geq 1$ let $\text{F}^n$ be the
$n$-fold composite $\text{F} \circ \text{F} \circ \cdots \circ
\text{F}$. We define $\text{F}_{n}' \colon \NIL(\mathcal{C};
\text{F}) \ra \NIL(\mathcal{C}; \text{F}^{n})$ to be the exact
functor which sends an object $c \colon C \ra \text{F}(C)$ to the
object given by the composite
$$
\xym@C=4em{ \text{F}^0(C)=C \ar[r]^-{\text{F}^0(c)=c} &
\text{F}(C) \ar[r]^-{\text{F}(c)} & \text{F}^2(C)
\ar[r]^-{\text{F}^2(c)} & \cdots \ar[r]^-{\text{F}^{n - 1}(c)} &
\text{F}^n (C) }
$$
and a morphism $f \colon (C,c) \ra (C',c')$ in
$\NIL(\mathcal{C};\text{F})$ to the morphism from $\text{F}_n'(C)$ to $\text{F}_n'(C')$ given by the underlying morphism $f
\colon C \ra C'$ in $\mathcal{C}$.
\item
 Let $U$ be a natural transformation $U \colon \text{F}^{\ell} \ra
\text{Id}_{\mathcal{A}}$ for some $\ell \in \nn$. For a natural
number $n$ we define 
$$
\text{F}_{\ell n + 1} \colon \NIL(\mathcal{C};\text{F}) \ra
\NIL(\mathcal{C};\text{F})
$$
to be the composition of the functor $\text{F}_{n \ell +
1}' \colon \NIL(\mathcal{C}; \text{F}) \ra
\NIL(\mathcal{C};\text{F}^{\ell n +1})$ and the functor
$\NIL(\text{Id},U^{n}) \colon \NIL(\mathcal{C};\text{F}^{n \ell
 +1}) \ra \NIL(\mathcal{C}; \text{F})$. The maps induced by
$\text{F}_{\ell n +1}$ on $\Nil_i(\mathcal{C};\text{F})$, for $i
\geq 0$, are also denoted by F$_{\ell n +1}$ and called
\emph{Frobenius operations}.
\end{enumerate}
\end{De}
\begin{Ex}[Frobenius on Farrell and Waldhausen Nil-groups]\label{ex}
Since $\al(x) = g x g^{-1}$
for some group element $g$, we can define an $RG$-module
homomorphism
$$
f \colon P \ot_{RG} (RG_{\al} \op X) \ra P
$$
by
$$
p \ot (r \op x) \mapsto p r g
$$
for $p \in P$, $r \in RG$ and $x \in X$.
The map $f$ induces a natural transformation $U^g$ between the functor
$\text{F}_{RG_{\al} \op X}$ and the identity.
\begin{enumerate}
\item  We obtain for every $n \in \nn$ and $i \in \nn$ a Frobenius operation on
the group $\Nil_i(RG;RG_{\al} \op X)$. 
\item As above we obtain a natural transformation $U^{g, g'}$ 
between the functor $\text{F}_{RG_{\al} \op X,RG_{\beta} \op Y}^2$ and the
identity. Thus we obtain for odd $n \in \nn$ and $i \in \nn$
a Frobenius operation on $\Nil_i(RG;RG_{\al} \op X,RG_{\beta} \op
Y)$. 
\item One can also define the Frobenius operations on various kinds of Waldhausen Nil-groups of generalized Laurent extension. Since we will not need this operations in the sequel we leaf the formulation of the precise statement to the reader.
\end{enumerate}
The Frobenius operations on lower Nil-groups are defined in
the obvious way.
\end{Ex}
\begin{Pro} \label{=01}
With the notation of the preceding definition, we obtain:
\begin{enumerate}
\item  For every $x \in \Nil_i(\mathcal{C};\text{F})$ there exists a natural number $L$ such that 
$$
\text{F}_n(x)= 0
$$
for $n \geq L$.
\item For $n_1$, $n_2 \in \nn$ whenever defined, we have
$$
\text{F}_{n_1} \text{F}_{n_2} = \text{F}_{n_1 \cdot n_2}.
$$
\end{enumerate}
\end{Pro}
\begin{proof}
For an exact category $\mathcal{E}$ we denote the category which is used in Quillen's
$Q$-construction is by $Q \mathcal{E}$. Let $Q_m \NIL(\mathcal{C};\text{F})$ be the full subcategories of the category $Q
\NIL(\mathcal{C};\text{F})$ consisting of
objects of nilpotency degree smaller or equal to $m$. For a
category $\mathcal{C}$ we denote the classifying space by $B
\mathcal{C}$. We have
\begin{align*}
K_i\big(\NIL(\mathcal{C};\text{F}) \big) & =
\pi_i \Omega B Q \NIL(\mathcal{C};\text{F}) \\
                           & = \colim_m \pi_i \Omega B Q_m
                           \NIL(\mathcal{C};\text{F})
\end{align*}
and
\begin{align*}
K_i ( \mathcal{C}) = \pi_i \Omega B Q_0
\NIL(\mathcal{C};\text{F}).
\end{align*}
We can find an $L$ such that $x \in \pi_i \Omega B Q_L \NIL(\mathcal{C};\text{F})$. Thus $\text{F}_n(x)=0$ for $n \geq L$.

The second identity follows since $U$ is a natural transformation and therefore for any map $c \colon C \rightarrow \text{F}(C)$ we have the following commutative diagram: 
\[
\xym@C=5em{\text{F}^{\ell  + 1}(C) \ar[r]^{U(\text{F}^{\ell+1}(C))} \ar[d]_{\text{F}^{\ell +1}(c^\ell)} & \text{F}(C) \ar[d]^-{\text{F}(c^\ell)} \\ \text{F}^{2\ell +1}(C) \ar[r]^{U(\text{F}^{2\ell +1}(C))} & \text{F}^{\ell+1}(C).}
\]
\end{proof}
\subsubsection{Verschiebung on Farrell Nil-Groups}
The main problem with the definition of Verschiebung operations is
that in general we do not have a map from $C$ to F$(C)$ which
plays the role of the identity. We use the assumption that we have
a natural transformation from the identity to F to get such a
map.
\begin{De}[Verschiebung]
Suppose that we have a natural transformation $U \colon
\text{Id}_{\mathcal{A}} \ra \text{F}$. Then we define, for $n \in \nn$,
$$
V_n \colon \NIL(\mathcal{C};\text{F}) \ra
\NIL(\mathcal{C};\text{F})
$$
by sending an object $c \colon C \ra \text{F}(C)$ to the object
$$
 \Big(C^{n},\text{%
\fontsize{8pt}{6pt}\selectfont $ \left(
\begin{array}{cccc} 0 & & & c
\\ U(C) & \ddots & &\\ & \ddots & 0 & \\ & & U(C) & 0 \end{array}
\right)$}\Big).
$$
A morphism $f$ from $(C,c)$ to $(C',c')$ is mapped to the morphism
$f^{\op n}$. The maps induced by these functors on
$\Nil_i(\mathcal{C};\text{F})$, for $i \geq 0$, are also denoted
by V$_{n}$ and called \emph{Verschiebung operations}.
\end{De}
To see that the functors V$_{n}$ are well-defined, note that for
every morphism $f$ from $C$ to $C'$ the diagram
$$
\xym{C \ar[r]^-{U(C)} \ar[d]_-f & \text{F} (C)
\ar[d]^-{\text{F}(f)} \\ C' \ar[r]^-{U(C')} & \text{F}(C')}
$$
commutes. In particular
$$
\xym@C=3em{C \ar[r]^-{U(C)} \ar[d]_-{c} & \text{F}(C)
\ar[d]^-{\text{F}(c)} \\ \text{F}(C) \ar[r]^-{U(\text{F}(C))} &
\text{F}^2(C)}
$$
commutes. This together with the fact that $c$ is nilpotent
implies that the morphism
$$
\text{%
\fontsize{8pt}{6pt}\selectfont $ \left(
\begin{array}{cccc}
0 & & & c \\
U(C) & \ddots & & \\
& \ddots & 0 & \\
& & U(C) & 0 \end{array} \right) $}
$$
is nilpotent.
\begin{Ex}[Verschiebung on Farrell Nil-groups]
Define the $RG$-module homomorphism
$$
f \colon P \ra P \ot_{RG} (RG_{\al} \op X)
$$
by
$$
p  \mapsto p \ot (g^{-1} \op 0)
$$
for $p \in P$. The map $f$ induces the required natural
transformation $U_g$ between the identity and $\text{F}_{RG_\al \op X}$.
Thus we obtain Verschiebung operations on $\Nil_i(RG;RG_\al \op
X)$ for $i \in \nn$. The Verschiebung operations on lower
Nil-groups are defined in the obvious way. Note that $U^g$ is
the left inverse of the natural transformation $U_g$ defined in Example \ref{ex}, i.e. the natural transformation $U^g \circ U_g$ is naturally isomorphic to the identity. If $X$ is the trivial module the
natural transformation $U^g$ is also a right inverse of $U_g$.
\end{Ex}
Consider $\nn$ with the multiplication. This gives $\nn$ the
structure of a semigroup. We define $\zz \nn$ to be the
''semigroup ring'' of $\nn$ with coefficients in $\zz$. The next
identity implies that we get a $\zz \nn$-module structure on
$\Nil_i(\mathcal{C};\text{F})$, where $n \in \nn$ operates on
$\Nil_i(\mathcal{C};\text{F})$ via $V_n$.
\begin{Pro} Let $n_1$ and $n_2$ be natural numbers.
With the assumptions of the preceding definition, we have
$$
V_{n_1} V_{n_2} = V_{n_1 \cdot n_2}
$$
as operations on $\Nil_i(\mathcal{C};\text{F})$.
\end{Pro}
\begin{proof}
This relation follows since there is a natural isomorphism between the
two functors. An application of Proposition
I. 1.3.1. in \cite{Wa2} yields the identity on
$K_i\big(\NIL(\mathcal{C};\text{F})\big)$ and therefore on
$\Nil_i(\mathcal{C};\text{F})$.
\end{proof}
\subsubsection{Verschiebung on Waldhausen Nil-Groups}
To define a Verschiebung operation on Waldhausen Nil-groups of
generalized free products we proceed in a similar manner as for
Farrell Nil-groups.
\begin{De}[Verschiebung on Waldhausen Nil-groups] \label{DefVer}
Let $\text{F}_1$ and $
\text{F}_2 \colon \mathcal{A} \ra \mathcal{A}$ be exact functors.
Let $\text{F}_W$ be the endofunctor of $\mathcal{A} \times
\mathcal{A}$ whose value at $(A,A')$ is
$(\text{F}_1(A'),\text{F}_2(A))$. Suppose that we have natural transformations
\begin{align*}
U_1 \colon \text{Id}_{\mathcal{A}} & \ra \text{F}_1 \\
U_2 \colon \text{Id}_{\mathcal{A}} & \ra \text{F}_2
\end{align*}
of exact functors $\mathcal{A} \ra \mathcal{A}$. For $\ell \in
\nn$ we define an exact functor
$$
V_{2 \ell + 1} \colon \NIL(\mathcal{C} \times \mathcal{C}
;\text{F}_W) \ra \NIL(\mathcal{C} \times \mathcal{C} ;\text{F}_W)
$$
which maps an object $(C,D,c,d)$ to
 \begin{center}
 \makebox[0pt]{
$(\text{%
\fontsize{7pt}{4pt}\selectfont $ (C \op D)^{\ell} \op C, $} \text{%
\fontsize{7pt}{4pt}\selectfont $(D \op C)^{\ell} \op D,$} \text{%
\fontsize{7pt}{4pt}\selectfont $ \left( \begin{array}{cccc}
 0 & & & c \\
 U_1(C) & \ddots & & \\
 & \ddots & 0 & \\
 & & U_1(D) & 0 \end{array} \right), $}
\text{%
\fontsize{7pt}{4pt}\selectfont $
 \left( \begin{array}{cccc}
 0 & & & d \\
 U_2(D) & \ddots & &\\
 & \ddots & 0 & \\
 & & U_2(C) & 0 \end{array}\right) $ }).
$}
\end{center}
If $(f,g)$ is a morphism from $(C,D,c,d)$ to $(C',D',c',d')$ of
$\NIL(\mathcal{C} \times \mathcal{C};\text{F}_W)$, we define V$_{2
\ell + 1}\big((f,g)\big)$ to be the morphism $\big((f \op g)^{\op
\ell} \op f, (g \op f)^{\op \ell} \op g \big)$.

The maps induced by these functors on $\Nil_i(\mathcal{C} \times
\mathcal{C} ;\text{F}_W)$, for $i \geq 0$, are also denoted by
V$_{2 \ell + 1}$ and called \emph{Verschiebung operations}.
\end{De}

The functors V$_{2 \ell + 1}$ are well-defined by the same
reasoning as above.
\begin{Ex}[Verschiebung on Waldhausen Nil-groups]
Since $U_g$ and $U_{g'}$ are natural transformations
between the identity and $\text{F}_{RG_{\al} \op X}$ and
$\text{F}_{RG_{\beta} \op Y}$ we obtain for every odd $n \in
\nn$ and $i \in \nn$ a Verschiebung operation on
$\Nil_i(RG;RG_{\al} \op X,RG_{\beta} \op Y)$. Verschiebung
operations on lower Nil-groups are defined in the obvious way.
Note that the natural transformation $U^{g,g'}$ 
is a left inverse of the natural transformation
induced by $U_g$ and $U_{g'}$. If $X$ and $Y$ are the trivial module
it is also a right inverse.
\end{Ex}
\begin{Pro}
Let the notation be as in the preceding definition. Let $n_1$ and
$n_2$ be odd natural numbers. With the assumptions of the
preceding definition, we have
$$
V_{n_1} V_{n_2} = V_{n_1 \cdot n_2}
$$
as operations on $\Nil_i(\mathcal{C} \times \mathcal{C};\text{F}_W)$.
\end{Pro}
\begin{proof}
The identity follows in the same manner as above.
\end{proof}
\subsection{Relations}\label{Relations}
In the first part of this section, we prove that the Frobenius and
Verschiebung operations satisfy the relations
$$
\text{F}_n \text{V}_n (x) = x \cdot n
$$
on $\Nil_i(\mathcal{C};\text{F})$ and a similar relation on
$\Nil_i(\mathcal{C} \times \mathcal{C};\text{F}_W)$. We use
this relation to prove that Nil-groups are
either trivial or not finitely generated as abelian groups.

In the second part the relation
$$
 \text{V}_n(y \ast \text{F}_n x) = (\text{V}_n y) \ast x.
$$
is proven.
\subsubsection{Non Finiteness Results}
\begin{De}[$\sigma$]
Let the notation be as in Definition \ref{DefVer}. We additionally
assume that $U_1$ and $U_2$ have a left inverse $U_1^{-1}$ and
$U_2^{-1}$. We define an exact endofunctor
\begin{align*}
\text{T} \colon \NIL(\mathcal{C} \times \mathcal{C}; \text{F}_W)
& \ra \NIL(\mathcal{C} \times \mathcal{C}; \text{F}_W) \\
(C,D,c,d) & \mapsto (D,C,U_1(C) \circ U_2^{-1}(C) \circ d,U_2(D)
\circ U_1^{-1}(D) \circ c)
\end{align*}
The induced map on $\Nil_i(\mathcal{C} \times \mathcal{C};
\text{F}_W)$ is denoted by $\sigma$.
\end{De}
Note that $\sigma$ is an group automorphism of order two.
\begin{Pro} \label{m.1} Let $\mathcal{A}$ be an abelian category and let $\mathcal{C} \subseteq \mathcal{A}$ be a full
subcategory which is closed under extension.
\begin{enumerate}
\item Let $\text{F} \colon \mathcal{A} \ra \mathcal{A}$ be an
exact functor. Suppose that we have
a natural transformation $U \colon \text{Id}_{\mathcal{A}} \ra
\text{F}$ which has a left inverse natural transformation. For $n \in
\nn$, we have
$$
\text{F}_n \text{V}_n (x) = x \cdot n
$$
for all $x$ in $\Nil_i(\mathcal{C};\text{F})$ and $i \in \nn$.
\item Let $\text{F}_1, \text{F}_2 \colon \mathcal{A} \ra
\mathcal{A}$ be exact functors. Suppose that we have natural transformations
$
U_1 \colon \text{Id}_{\mathcal{A}}  \ra \text{F}_1 
$ and $
U_2 \colon \text{Id}_{\mathcal{A}}  \ra \text{F}_2
$
which have left inverse natural transformations $U_1'$ and $U_2'$ respectively. The natural transformations $U_1'$ and $U_2'$ induce a natural transformation between $\text{F}_W^2$ and the identity and therefore we have Frobenius and Verschiebung operations on $\Nil_i(\mathcal{C} \times \mathcal{C};\text{F}_W)$. For odd $n=2 \ell + 1$, we have
$$
\text{F}_n \text{V}_n (x) = x \cdot (\ell +1) + \sigma (x) \cdot
\ell
$$
for all $x$ in $\Nil_i(\mathcal{C} \times \mathcal{C};\text{F}_W)$
and $i \in \nn$.
\end{enumerate}
\end{Pro}
\begin{proof}
For $\NIL(\mathcal{C};\text{F})$ the natural transformation $U$ induce a natural transformation between F$_n \text{V}_n$ and the functor which sends an object to the $n$-th fold direct sum. The induce natural transformation is an isomorphism if $U$ has a left inverse.

Similar arguments work for $\NIL(\mathcal{C} \times \mathcal{C};\text{F}_W)$.
\end{proof}
\begin{Le}
Let $G$ be an abelian group with finite torsion subgroup $T$. Let
$\sigma$ be a group automorphism of order two. If $n \in \nn$ is
$-1$ modulo $|T|$ then
$$
\Phi_{n}(x) := x \cdot (n + 1) + \sigma(x) \cdot n
$$
is a monomorphism of $G$ into itself.
\end{Le}
\begin{proof} First we prove that $\Phi_n(x)=0$ implies $x \in
T$. We define $u_x:= \sigma (x) + x$ and $v_x(x):= - \sigma (x)
+x$. We have $u_x + v_x  = 2 x$, $\sigma(u_x)  = u_x$ and $\sigma(v_x)  = - v_x.$
Since $\Phi_n(x)=0$ implies that $\Phi_n(2 x)=0$, we have
\[
0  = \Phi_n(u_x + v_x)  = (2 n +1) u_x + v_x.
\]
Applying $\sigma$ to this equality yields $(2 n +1) u_x -v_x  = 0.$
Hence $u_x \in T$. Therefore, $v_x \in T$ and also $u_x + v_x \in
T$. In particular $2 x \in T$ and $x \in T$.

Since $|T|$ divides $n+1$ we have
\[
\Phi_n(x) = - \sigma(x) \,\,\text{for} \,\, x \in T.
\]
Therefore $\Phi_n$ is a monomorphism on $T$. This implies the lemma since, by the first part of this proof, all elements in the kernel of $\Phi_n$ are in $T$.
\end{proof}
\begin{Co} Let the notation be as in the preceding proposition.
\label{infinite} The groups $\Nil_i(\mathcal{C};\text{F})$
and~$\Nil_i(\mathcal{C} \times \mathcal{C};\text{F}_W)$ are either
trivial or not finitely generated as abelian groups for $i \in
\nn$.
\end{Co}
\begin{proof}
We prove the result for $\Nil_i(\mathcal{C} \times
\mathcal{C};\text{F}_W)$, similar arguments work for
$\Nil_i(\mathcal{C};\text{F})$.

If $\Nil_i(\mathcal{C} \times \mathcal{C};\text{F}_W)$ is finitely generated, its torsion subgroup $T$ is finite. By Proposition \ref{=01}, $\text{F}_n=0$ on $\Nil_i(\mathcal{C} \times \mathcal{C};\text{F}_W)$ for large $n$. Choosing $n$ large and congruent to $-1$ modulo $|T|$, $\text{F}_n \text{V}_n$ is a monomorphism by Proposition~\ref{m.1}, contradiction.   
\end{proof}
\begin{Co}
Let $R$  be a ring, let $G$ be a group, let $X$ and $Y$ be
arbitrary $RG$-bimodules and let $\al$ and $\beta$ be inner group
automorphisms of $G$. The groups $\Nil_i(RG;RG_\al \op X)$ and
$\Nil_i(RG;RG_{\al} \op X,RG_{\beta} \op Y)$ are either trivial or
not finitely generated as abelian groups for $i \in \zz$.
\end{Co}
\subsubsection{The Relation V$_n(y \ast $F$_n x) = ($V$_n y) \ast x$
on Nil-Groups}
The proof of the relation
$$
\text{V}_n(y \ast \text{F}_n x) = (\text{V}_n y) \ast x
$$
requires a little bit more machinery. The basic idea, which is due
to Stienstra \cite{Sti}, is to define some exact category
$\END(\Delta;S_n)$ and an exact pairing
\begin{align*}
\END(\Delta;S_n) \times \END(R) \times \NIL(\La;X,Y) & \lra
\NIL(\La;X,Y)\\
c \times x \times y & \lra (c,x,y).
\end{align*}
Then we prove that there are objects $C_1$ and $C_2$ of
$\END(\Delta;S_n)$ such that V$_n(x \times \text{F}_n y) = (C_1,
x,y)$ and $(\text{V}_n x) \times y =(C_2,x,y)$. Using the product
map developed by Waldhausen \cite{Wa1,Wa4}, we can now prove the
relation by proving that $[C_1]=[C_2]$ as elements in
$\End_0(\Delta;S_n)$. We will prove the identity just for Waldhausen Nil-groups a similar pairing can be defined for Farrell Nil-groups. 

To start this program, we define the category $\END(R;T)$.
\begin{De}[$\END(R;T)$]
Let $R$ be a ring and let $T \subset R[t]$ be a multiplicatively
closed set containing $t$. Recall that $\text{F}_{R,R}$ is the functor from $\Mod(R) \times \Mod(R)$ to $\Mod(R) \times \Mod(R)$ whose value at the pair $(M,N)$ is $(N,M)$. We define
$\END(R;T)$ to be the full exact subcategory of
$\END(\mathbf{P}(R) \times \mathbf{P}(R);\text{F}_{R,R})$ consisting of quadruples $(P,Q,p,q)$ such that there
exist polynomials $p_1(t)=\sum t^i \la_i $ and $p_2(t)=\sum t^i
\mu_i $ in $T$ with $\sum (p \circ q)^i \la_i = 0$ and $\sum (q
\circ p)^i \mu_i =0$.
\end{De}

Let $\delta$ be the ring automorphism of the polynomial ring $\zz[y,v,w]$ induced by
mapping $v$ to $w$ and $w$ to $v$. Let $S_n$ be the
multiplicatively closed subset of the
twisted polynomial ring  $\Delta :=\zz [y,v,w]_{\delta}[x]$ generated by
$t$, $t^n - v^{2n} \cdot y^2 \cdot x^{2n}$ and $t^n - w^{2 n}
\cdot y^2 \cdot x^{2n}$. 

For odd $n$ objects of $\END(\Delta;S_n)$ are for example 
\begin{center}
\makebox[0pt]{%
$C_1 =(\text{\fontsize{7pt}{5pt}\selectfont $ \Delta^n, \Delta^n,$}
\text{%
\fontsize{7pt}{5pt}\selectfont $ \left(
\begin{array}{ccccc}
0 & & & & y x^n v\\
v & 0 & & & \\
  & w & \ddots & & \\
  &   & \ddots & 0 &  \\
  &   &        & w & 0
\end{array} \right)$},
\text{%
\fontsize{7pt}{5pt}\selectfont $ \left(
 \begin{array}{ccccc}
0  & & & & y x^n w \\
w  & 0 & & & \\
   & v & \ddots & & \\
   &   & \ddots & 0 & \\
   & & & v & 0
\end{array} \right)$})
$}
\end{center}
and
\begin{center}
\makebox[0pt]{
$C_2 = (\text{\fontsize{7pt}{5pt}\selectfont $ \Delta^n,\Delta^n,$} \text{%
\fontsize{7pt}{5pt}\selectfont $ \left(
\begin{array}{ccccc}
  0 & & & &  y x v \\
   x v & 0 & & & \\
    & x v & \ddots & & \\
    & & \ddots & 0 & \\
    & & & x v & 0
 \end{array} \right),$} \text{%
\fontsize{7pt}{5pt}\selectfont $ \left(
\begin{array}{ccccc}
 0 & & & &  y x w \\
  x w & 0 & & &\\
   &  x w & \ddots & & \\
   & & \ddots & 0 & \\
   & & &  x w & 0
 \end{array} \right)$} ).$}
\end{center}
which are annihilated by $(t^n - v^{2n} \cdot y^2 \cdot x^{2n})
\cdot (t^n - w^{2 n} \cdot y^2 \cdot x^{2n})$. This means the following. Put $t$ equal to the product of the two morphism of either $C_1$ or $C_2$. The polynomial $(t^n - v^{2n} \cdot y^2 \cdot x^{2n})
\cdot (t^n - w^{2 n} \cdot y^2 \cdot x^{2n})$ vanishes. To confirm this note that $t^n$ is obviously a diagonal matrix. A little computation shows that $t^n$ is in fact a diagonal matrix with entries given by the above values. 

Let $\La$ be an algebra over a ring $R$, let $X$ and $Y$ be left
flat $\La$-bimodules and let $U_X$ and $U_Y$ be natural
transformations between the identity and $\text{F}_{X}$ and
$\text{F}_Y$ with right inverse natural transformation $U_X^{-1}$
and $U_Y^{-1}$. We define an exact pairing
\begin{center}
\makebox[0pt]{ $\END(\Delta;S_n) \times \END(R)
\times \NIL(\La;X,Y) \ra \NIL(\La;X,Y) $}
\end{center}
where
$$
(C,D,c,d) \times (B,\varphi) \times (P,Q,p,q)
$$
is mapped to the object
\begin{center}
\makebox[0pt]{ $(\text{%
\fontsize{9pt}{6pt}\selectfont $ C \ot_{\Delta} B \ot_R (P \op
Q),D \ot_{\Delta} B \ot_R (P \op Q),
c \ot \text{id} \ot U_X(P \op Q),d \ot \text{id} \ot U_Y(P \op Q)$}),$}
\end{center}
where $x$ acts on $B \ot_{R} (P \op Q)$ via $\text{id} \ot \left(
\begin{array}{cc} 0 & U_Y^{-1}(Q) \circ q \\ U_X^{-1}(P)
\circ p & 0 \end{array} \right) $, $y$ acts via $\varphi \ot
\text{id}$, $v$ acts via $\text{id} \ot \left( \begin{array}{cc}
 1 & 0
\\ 0 & 0 \end{array} \right)$ and $w$ acts via $\text{id} \ot \left(
\begin{array}{cc} 0 & 0
\\ 0 & 1 \end{array} \right)$.
\begin{De}[$\END^{D1}(\La;T)$]
Let $\La$ be a ring and let $T \subseteq \La[t]$ be a
multiplicatively closed set containing $t$. Let
$\END^{D1}(\La;T)$ be the full subcategory of
$\END(\Mod(\La) \times \Mod(\La);\text{F}_{\La,\La})$ consisting of those objects
which have an $\END(\La;T)$-resolution of
length one.
\end{De}
In the same way that the Lemma \ref{7} is proven we obtain that $\END^{D1}(\La;T)$ is an exact category.
\begin{De}[$\End_0(\La;T)$,$\End_0^{D1}(\La;T)$]
Let the notation be as above.
\begin{enumerate}
\item We define $\End_0(\La;T)$ to be the
kernel of the map on $K_0$ which is induced by the forgetful
functor from $\END(\La;T)$ onto
$\mathbf{P}(\La) \times \mathbf{P}(\La)$ whose value at
$(P,Q,p,q)$ is $(P,Q)$.
\item Let $\mathbf{P}^{D1}(\La)$ be the full subcategory of
$\Mod(\La)$ consisting of those objects which have a
$\mathbf{P}(\La)$-resolution of length one.
\item We define $\End_0^{D1}(\La;T)$ to be the
kernel of the map on $K_0$ which is induced by the forgetful
functor from $\END^{D1}(\La;T)$ onto
$\mathbf{P}^{D1}(\La) \times \mathbf{P}^{D1}(\La)$ whose value
at $(M,N,m,n)$ is $(M,N)$.
\end{enumerate}
\end{De}
\begin{Le} \label{c1=c2}
Let the notation be as above. The inclusion map
$$
\End_0(\La;T) \lra \End_0^{D1}(\La;T)
$$
is an isomorphism.
\end{Le}
\begin{proof}
The Resolution Theorem
gives that the inclusion maps from
$K_0\big(\mathbf{P}(\La)\big)$ into
$K_0\big(\mathbf{P}^{D1}(\La)\big)$ and from
$K_0\big(\END(\La;T)\big)$ into
$K_0\big(\END^{D1}(\La;T)\big)$ are isomorphisms.
Combining these two results we get the lemma.
\end{proof}
\begin{Le}
We have
\[
[C_1] = [C_2]
\] 
in $\End_0(\Delta;S_n)$.
\end{Le}
\begin{proof}
 Since the groups
$\End^{D1}_0(\Delta;S_n)$ and
$\End_0(\Delta;S_n)$ are
naturally isomorphic it is enough to prove that the objects
represent the same element in $\End^{D1}_0(\Delta;S_n)$.

Consider the map $\iota \colon C_1 \ra C_2$ in
$\End^{D1}_0(\Delta;S_n)$
which is induce by the injective maps
\[
\iota_1 := \text{%
\fontsize{8pt}{6pt}\selectfont $ \left(
\begin{array}{cccc}
1  & & & \\
& x & & \\
& & \ddots & \\
 & & & x^{n-1}
\end{array} \right)$} \,\, \text{and} \,\,
\iota_2 := \text{%
\fontsize{8pt}{6pt}\selectfont $ \left(
\begin{array}{cccc}
 1  & & & \\
& x & \\
& & \ddots & \\
 & & & x^{n-1}
\end{array} \right)$}.
\]

The object
\begin{center}
\makebox[0pt]{
$( \text{%
\fontsize{10pt}{8pt}\selectfont $   \bigoplus_{i=0}^{n-1} \Delta/x^i,\bigoplus_{i=0}^{n-1} \Delta/x^i,$}
\text{%
\fontsize{6pt}{4pt}\selectfont $ \left(
\begin{array}{ccccc}
  0 & & & &  \\
    & 0 & & & \\
    & x v & \ddots & & \\
    & & \ddots & 0 & \\
    & & & x v & 0
 \end{array} \right),$} \text{%
\fontsize{6pt}{4pt}\selectfont $ \left(
\begin{array}{ccccc}
 0 & & & &  \\
   & 0 & & &\\
   & x w & \ddots & & \\
   & & \ddots & 0 & \\
   & & & x w & 0
 \end{array} \right)$})
$}
\end{center}
 is the cokernel of $\iota$ in $\END^{D1}(\Delta;S_n)$. In the following we denote this object by $C$. We get a short exact sequence
$$
\xym{0 \ar[r] & C_1 \ar[r]^-{\iota} & C_2 \ar[r] & C \ar[r] & 0}
$$
in $\END^{D1}(\Delta;S_n)$. In $\End_0^{D1}(\Delta;S_n)$ objects of the form 
\[
\big[ \big( \bigoplus_{i=0}^n
M_i,\bigoplus_{i=0}^n M_i,m_1,m_2\big)\big]
\]
 with $m_1$ and $m_2$
lower triangular matrices vanish. Thus
$
\big[ \big( C \big) \big]=0
$
in the group $\End_0^{D1}(\Delta;S_n)$, which implies that
$[C_1 ]$ and $[ C_2 ]$ represent the same element in
$\End_0^{D1}(\Delta;S_n)$ and
therefore the required identity.  
\end{proof}
\begin{Th}\label{vm=fm} Let $\La$ be an algebra over a ring $R$,
let $X$ and $Y$ be left flat $\La$-bimodules and let $U_X$ and
$U_Y$ be natural transformations between the identity and
$\text{F}_{X}$ and $\text{F}_Y$ with right inverse natural
transformations $U^{-1}_X$ and $U^{-1}_Y$. Let $y$ be an element
of $\End_0(R)$ and let $x$ be an element of $\Nil_i(\La;X,Y)$ for
$i \in \zz$. For odd $n \in \nn$, we have
$$
\text{V}_n(y \ast \text{F}_n x) = (\text{V}_n y) \ast x
$$
where $\ast$ is the module multiplication defined in
Section~\ref{pairing}.
\end{Th}
\begin{proof} We prove the statement for higher Nil-groups. The arguments extend
to lower Nil-groups in the obvious way. The pairing between the
categories $\END(R) \times \NIL(\La;X,Y)$ and $\END(\Delta;S_n)$ defined above gives
that every object $(C,D,c,d)$ of
$\END(\Delta;S_n)$ defines an exact functor from $\END(R) \times \NIL(\La;X,Y)$ to $\NIL(\La;X;Y)$. Since objects
of the form $(C,D,0,0)$ reflect the category $\END(R) \times \NIL(\La;X,Y)$
into the subcategory where elements are of the form $(P,Q,0,0)$ we
get that every element in
$\End_0(\Delta;S_n)$ gives
rise to a map from $\End_0(R) \times \Nil_i(\La;X,Y)$ to $\Nil_i(\La;X,Y)$. Let $(B,\varphi)$ be an object of $\END(R)$, let $(P,Q,p,q)$ be an
object of $\NIL(\La;X,Y)$ and let $n=2 \ell +1$ for $\ell \in \nn$.
We have
\begin{align*}
& \text{V}_n \big((B,\varphi) \times \text{F}_n (P,Q,p,q) \big)
\op (B \ot  Q,B \ot  P,0,0)^{\ell +1} \op (B \ot P,B \ot
Q,0,0)^{\ell} \\ \\
& \cong  C_1 \times (B,\varphi) \times (P,Q,p,q)
\end{align*}
and
\[
 \big(\text{V}_n (B,\varphi) \big) \times (P,Q,p,q) \op (B\ot Q,B
\ot  P,0,0)^n  \cong C_2  \times (B,\varphi) \times (P,Q,p,q) .
\]
Thus it remains to show that $C_1$ and $C_2$ represent the same elements in the group $\End_0(\Delta;S_n)$. But this is the statement of Lemma \ref{c1=c2}. 
\end{proof}
\subsection{Nil-Groups as Modules over the Ring of Witt Vectors}
As a direct corollary of the relation V$_n(y \ast $F$_n x) = ($V$_n y) \ast x$  we obtain that Nil-groups are modules over the ring of Witt vectors. Let us briefly recall the
definition of the ring of Witt vectors. For a commutative ring $R$ the \emph{ring of (big) Witt
vectors} is the ring $1 + tR \llbracket t \rrbracket$ of power
series with constant term $1$. The underlying additive group of
the ring of Witt vectors is the multiplicative group of $1 + t R
\llbracket t \rrbracket$. The multiplication is the unique
continuous functorial operation $*$ for which
$$
(1-at) \ast (1-bt) = (1-abt)
$$
holds for all $a$, $b \in R$. In the sequel, the ring of Witt
vectors is denoted by $W(R)$. We define ideals~$I_N : = (1 + t^N R
\llbracket t \rrbracket)$ for all $N \in \nn$. The resulting
topology on the ring of Witt vectors is called the \emph{$t$-adic
topology}. Let $W_N(R) := W(R) / I_N$ be the \emph{truncated ring of Witt
vectors}.

Since $R$ is assumed to be
commutative, the tensor product induces a ring structure on
$\End_0(R)$. The characteristic polynomial
$$
\chi \big((B,\varphi)\big) : = \det (\text{id}_B - t \cdot
\varphi)
$$
defines a map from $\End_0(R)$ into $W(R)$. 
A theorem of Almkvist states that this map is an
injective ring homomorphism whose image is dense in the ring of
Witt vectors with respect to the $t$-adic topology \cite{Al2}. This result implies that to extend the $\End_0(R)$-module structure to a $W(R)$-module structure we need to show that for every $x \in \Nil$ there exist an arbitrary $N \in \nn$ such that $x$ is annihilated by $\End_0(R) \cap I_N$. Since $I_M$ is contained in $I_N$ for $M \geq N$
we also obtain that $\End_0(R) \cap I_M$ annihilates $x$ for $M \geq N$.
\begin{Th}\label{wittvec}
Let $\La$ be an algebra over a commutative ring $R$, let $X$ and
$Y$ be left flat $\La$-bimodules. If we have natural
transformations $U_X$ and $U_Y$ between the identity and
$\text{F}_{X}$ and $\text{F}_Y$ which have right inverse natural
transformations then for every x in $\Nil_i(\La;X)$ or
$\Nil_i(\La;X,Y)$ there is an N such that x is annihilated by $\End_0(R) \cap I_N$. Thus $\Nil_i(\La;X)$ and $\Nil_i(\La;X,Y)$ are
modules over the ring of Witt vectors.

Moreover this Witt vector module structure is continuous, meaning that for every finitely generated submodule $M$ we can find an $N$ such that the $W(R)$-module structure restricts to a $W_N(R)$-module structure.
\end{Th}
\begin{proof}
We prove the results for Waldhausen Nil-groups similar arguments work for the Farrell Nil-groups.

Let $x$ be an element of $\Nil_i(\La;X,Y)$. Proposition~\ref{=01} implies that there is an odd $N \in \nn$ such that $\text{F}_N(x) = 0$. Let $y$ be an element of $\End_0(R) \cap I_N$. Since the Verschiebung operation on $W(R)$ restrict to the Verschiebung operation on $\End_0(R)$ \cite{Gr2} we have $y = \text{V}_N(y')$ for some $y' \in \End_0(R)$. Theorem~\ref{vm=fm} implies now
\[
y \ast x = (\text{V}_N y') \ast x = \text{V}_N (y' \ast \text{F}_N x) = 0
\]   
and therefore $(\End_0(R) \cap I_N) \ast x =0$.

The second part of the theorem follows in the same way.
 \end{proof}
\section{Applications} 
We combine the results of the preceding sections to show that taking Nil-groups and localization commutes. The main application is torsion results. As a consequence, of the torsion results, we obtain that the relative assembly map from the family of finite subgroups to the family of virtually cyclic subgroups is rationally an isomorphism.
\subsection{Nil and Localization Commute}
\label{Imlocalization}
\begin{Le}\label{Wittloca}
Let $\La$ be an algebra over a commutative ring $R$ and let $X$
and $Y$ be $\La$-bimodules such that there are natural
transformations from the identity to the functor F$_X$ and F$_Y$
which have right inverse natural transformations. Then for every
multiplicatively closed set $S \subset R$ of non zero divisors
satisfying $s \cdot x = x \cdot s$ for all $s \in S$ and $x \in X$ (or $x \in Y$) we have
$$
\{(1-st) \, | \, s \in S \}^{-1} W(R) \ot_{W(R)} \Nil_i(\La;X)
\cong \Nil_i(\La_S;{}_S X_S)
$$
and
$$
\{(1-st) \, | \, s \in S \}^{-1} W(R) \ot_{W(R)} \Nil_i(\La;X,Y)
\cong \Nil_i(\La_S;{}_S X_S,{}_S Y_S),
$$
for all $i \in \zz$.
\end{Le}
\begin{proof}
Suppose first that $S$ is generated by one element $s$, so that multiplication by $(1-st)$ is induced by the functor S defined after Definition \ref{Imd1}. In this case the result is just a restatement of Corollary \ref{Th}. Iterating yields the result when $S$ is finitely generated; for a general multiplicatively closed $S$, we take the colimit.
\end{proof}
\begin{Th}
\label{localization}Let $R$ be $\zz_T$ for some multiplicatively
closed set $T \subseteq \zz-\{0\}$, $\hat{\zz}_p$ or a commutative
$\qq$-algebra and let $\La$ be an $R$-algebra. Let $X$ and $Y$ be
$\La$-bimodules such that there are natural transformations from
the identity to the functor F$_X$ and F$_Y$ which have a right
inverse. Then for every multiplicatively closed set $S \subset R$
of non zero divisors satisfying $s \cdot x = x \cdot s$ for all $s
\in S$ and $x \in X$ (or $x \in Y$) there are
isomorphisms of $R_S$-modules
$$
R_S \ot_{R} \Nil_i(\La;X) \cong \Nil_i(\La_S;{}_S X_S)
$$
and
$$
R_S \ot_{R} \Nil_i(\La;X,Y) \cong \Nil_i(\La_S;{}_S X_S,{}_S Y_S),
$$
for all $i \in \zz$.
\end{Th}
\begin{proof} This proof follows Weibel's proof of the same
identity for Bass Nil-groups \cite[page 489]{We1}. Again just the
groups $\Nil_i(\La;X)$ are treated, the same arguments work for
$\Nil_i(\La;X,Y)$.

The group $\Nil_i(\La;X)$ is a colimit over the family of finitely
generated $W(R)$-submodules~$M$. Since $M$ is assumed to be
finitely generated, Theorem~\ref{wittvec} 
implies that the $W(R)$-module structure restricts to a
$W_N(R)$-module structure for a certain $N$. For a ring $R$ which
is $\zz_T$ for some multiplicatively closed set $T \subseteq
\zz-\{0\}$, $\hat{\zz}_p$ or a commutative $\qq$-algebra Weibel
\cite[Proposition 6.2]{We1} proves the identity
$$
\{(1-st) \, | \, s \in S \}^{-1} W_N(R) \cong W_N(R_S).
$$
Since $R$ is a $\la$-ring, $M$ caries an $R$-module structure. We
have
\begin{align*}
 R_S \ot_R M  & \cong  W_N(R_S) \ot_{W(R)} M \\
 & \cong \{(1-st) \, | \, s \in S \}^{-1} W_N(R)
\ot_{W(R)} M \\
 & \cong \{(1-st) \, | \, s \in S \}^{-1} W(R) \ot_{W(R)} M.
\end{align*}
Taking the colimit of both sides gives
$$
 R_S \ot_R \Nil_i(\La;X) \cong \{(1-st) \, | \, s \in S \}^{-1} W(R)
 \ot_{W(R)} \Nil_i(\La;X).
$$
Thus Lemma~\ref{Wittloca} implies
$$
R_S \ot_R \Nil_i(\La;X) \cong \Nil_i(\La_S;{}_S X_S). \qedhere
$$
\end{proof}
\begin{Rem}
Note that the statement stays valid if the ring $R$ is a $\la$-ring and satisfies $\{(1-st) \, | \, s \in S \}^{-1} W_N(R) \cong W_N(R_S)$.
\end{Rem}
\subsection{Induction and Transfer on Nil-groups}\label{Torsionresults}
To prove the torsion results stated in the introduction we need to define induction and transfer maps on the Nil-groups.
\begin{De}[Induction and Transfer]
Let $\Gamma$ be a ring with a subring $\La$. Let $\iota \colon \La
\hookrightarrow \Gamma$ be the inclusion map and let $\al$ and
$\beta$ be a ring automorphisms of $\Gamma$ which restricts to
$\La$.
\begin{enumerate}
\item Define a functor $u$ from $\Mod(\La)$ to
$\Mod(\Gamma)$ which sends a $\La$-module $M$ to the
$\Gamma$-module $M \ot_{\La} \Gamma$. We denote the natural
transformation between $u \circ \text{F}_{\La_\al}$ and
$\text{F}_{\Gamma_\al} \circ u$ which is induced by the map
\begin{align*}
 \La_{\al} \ot_{\La} \Gamma & \ra \Gamma_{\al}\\
 \la \ot \gamma & \mapsto \la \cdot \al(\gamma),
\end{align*}
by $U$. We define $\iota_{i}$ to be $\Nil_i(u,U) \colon
\Nil_i(\La;\al) \ra \Nil_i(\Gamma;\al)$. 
\item We denote the natural transformation between $u \circ
\text{F}_{\La_\al,\La_\beta}$ and
$\text{F}_{\Gamma_\al,\Gamma_\beta} \circ u$ which is induced by
the maps
\begin{align*}
 \La_{\al} \ot_{\La} \Gamma & \ra \Gamma_{\al} \\
 \la \ot \gamma & \mapsto \la \cdot \al(\gamma),\\
 \La \ot_{\La} \Gamma & \ra \Gamma \\
 \la \ot \gamma & \mapsto \la \cdot \beta (\gamma),
\end{align*}
by $U_W$. We define $\iota_{i}$ to be $\Nil_i(u,U_W) \colon
\Nil_i(\La;\La_\al,\La_\beta) \ra
\Nil_i(\Gamma;\Gamma_\al,\Gamma_\beta)$. 
\end{enumerate}
We call these maps
\emph{induction maps}.

If we have additionally that $\Gamma_{\iota}$ is a finitely
generated projective right $\La$-module, we can define transfer
maps.
\begin{enumerate}
\item We define a functor T from $\NIL(\Gamma;\al)$ to
$\NIL(\La;\al)$ whose value at an object $(P,\nu)$ is
$(P_{\iota},\nu)$. On morphisms, T is the identity. 
\item We define a functor T from
$\NIL(\Gamma;\Gamma_{\al},\Gamma_{\beta})$ to
$\NIL(\La;\La_{\al},\La_{\beta})$ whose value at an
object~$(P,Q,p,q)$ is $(P_{\iota},Q_{\iota},p,q)$. On morphisms, T
is the identity. 
\end{enumerate}
The maps which are induced on the $i$-th Nil-groups is denoted by $\iota^i$ and
called \emph{transfer map}.
\end{De}
\begin{Ex}
The situation which will become important is that we have a group
$G$ and a group automorphism $\al$ of finite order $m$. In this
case we can form the semidirect product of $G$ and the cyclic
group of order $m$, which is denoted by $C_m$. Conjugation with
the element $(\text{id},1)$ extends $\al$ to an inner group
automorphism of $G \rtimes_{\al} C_m$, which is denoted by
$\tilde{\al}$. The group ring $R G \rtimes_{\al} C_m$, seen as a
bimodule over $R G$, is isomorphic to $\bigoplus_{i=0}^{m-1} R
G_{\al^i}$. Thus we get induction maps
$$
\iota_i \colon \Nil_i(RG;\al) \ra \Nil_i(R G \rtimes
C_m;\tilde{\al})
$$
and transfer maps
$$
\iota^i \colon \Nil_i(R G \rtimes C_m;\tilde{\al})  \ra
\Nil_i(RG;\al).
$$

For the definition of transfer and induction maps on Waldhausen Nil-groups of a generalized free product we will restrict to the case $\Nil_i(RG;RG_{\al},RG)$. We can do so without lost of generality since we have 
\[
\Nil_i(\Lambda;\Lambda_{\al},\Lambda_\beta) \cong \Nil_i(\Lambda;\Lambda_{\al \beta},\Lambda).
\]
for arbitrary ring automorphisms $\al$ and $\beta$ \cite[Proposition 3.2]{CP}.
By the same reasoning as above we can define induction maps
$$
\iota_i \colon \Nil_i(RG;RG_{\al},RG)  \ra \Nil_i(RG
\rtimes C_m;RG \rtimes C_{m \, \tilde{\al}},
RG \rtimes C_{m})
$$
and transfer maps
$$
\iota^i \colon \Nil_i(RG \rtimes C_{m};RG \rtimes C_{m \, \tilde{\al}}, RG \rtimes C_{m}) \ra \Nil_i(RG;RG_{\al},RG).
$$
\end{Ex}
\begin{Le} \label{Transfer}
Let the notation be as in the preceding example. We have
$$
\iota^i \circ \iota_i (x) =x \cdot m,
$$
for $x$ in $\Nil_i(R G;\al)$ or $\Nil_i(RG;R G_{\al},R G)$ and $i \in \zz$. 
\end{Le}
\begin{proof}
We prove the lemma for $\Nil_i(R G;\al)$ similar arguments work for the group $\Nil_i(RG;R G_{\al},R G)$.

We have 
\[
\text{T} \circ \NIL(U,u) \big((P,p)\big) = (\op_{i=0}^{m-1} P_{\al^i}, \op_{i=0}^{m-1} p ). 
\]
Define a functor
\begin{align*}
\text{G} \colon \NIL (\mathbf{P}(\zz G); \zz G_\al) & \ra \NIL(\mathbf{P}( \zz G); \zz G_\al) \\
(P,p) & \mapsto (P_{\al},p).
\end{align*}
In the sequel we prove that $\text{G}$ induces the identity on $\Nil_i(\zz G; \zz G_\al)$. Iterated use of this argument proves that $\text{T} \circ \NIL(U,u)$ induces multiplication by $m$ and therefore the lemma.

Define functors 
\begin{align*} 
\text{I} \colon \NIL (\mathbf{P}(\zz G); \zz G_\al) & \ra \NIL(\mathbf{P}( \zz G); \zz G_\al) \\
(P,p) & \mapsto (P_{\al^i} \op P_{\al^{i+1}},\left(
\begin{array}{cc}
  p & \text{Id} \\
  0  & 0 
 \end{array} \right))\\ 
\text{J} \colon \NIL (\mathbf{P}( \zz G); \zz G_\al ) & \ra \NIL(\mathbf{P}( \zz G); \zz G_\al ) \\
(P,p) & \mapsto (P_\al,0)\\ 
\Ker \colon \NIL (\mathbf{P}(\zz G); \zz G_\al) & \ra \NIL(\mathbf{P}( \zz G); \zz G_\al ) \\
(P,p) & \mapsto \big(\Ker\big(\left(
\begin{array}{cc}
  p & \text{Id} 
 \end{array} \right)\big),0). 
\end{align*}

We have two exact sequences of exact functors: 
\[
\xym{0 \ar[r] & \text{Id} \ar[r] & \text{I} \ar[r] & \text{J} \ar[r] & 0}
\]
\[
\xym{0 \ar[r] & \Ker \ar[r] &  \text{I} \ar[r] & \text{G} \ar[r] & 0.}
\]
The Additivity Theorem together with the observation that J and $\Ker$ induce the trivial map on $\Nil_i(RG;RG_\al)$ implies that Id and G induce the same map on $\Nil_i(RG;RG_\al)$ and therefore the lemma. 
\end{proof}
\subsection{Torsion Results}\label{Torsionresults}
In view of the result that Nil and localization commute it is important to find a class of groups such that $\zz[1/n] G$ is regular. In the first part of this section we prove that for every polycyclic-by-finite group we can find such an $n$. In the second part we apply this results to obtain torsion results for the Nil-groups of such groups.
\begin{Le}
\label{fi} Let $R$ be a ring, let $G$ be a group and let $H$ be a
subgroup of finite index $n$ with the property that $R H$ is
regular. If n is invertible in $R$, then $R G$ is also regular.
\end{Le}
\begin{proof}
 The ring $R G$ is right noetherian since $R H$ is right noetherian and $R
 G$ is a finitely generated right module over $R H$. Let $M$ be a
 finitely generated $R G$-module. The module $M$ seen as an $R H$-module
 is denoted by $\text{res}\,M$.
 Since $R H$ is regular, $\text{res}\,M$ has a finite projective $R H$-resolution.
 Applying $- \ot_{R H} R G $ yields an $R G$-resolution of
 $\text{res}\, M \ot_{R H} R G$. Let $S$ be a set of representatives of right
 cosets. Since $n$ is invertible in $R$, we can define the following
 map:
\begin{align*}
 M & \ra  \text{res}\, M \ot_{R H} R G \\
 m & \mapsto 1/n \sum_{g \in S}  m  g^{-1} \ot g.
\end{align*}
The prove that this maps is an an $R G$-module map is left to the reader.
Since the composition of this map with the canonical map from
$\text{res}\, M \ot_{R H} R G$ to $M$ is the identity, we get that
$M$, as an $R G$-module, is a direct summand of $\text{res}\, M
\ot_{R H} RG$. Hence $M$ has a finite projective resolution.
\end{proof}
Every polycyclic-by-finite group $G$ contains poly-infinite cyclic subgroup of finite index \cite[1.5.12]{McC}.
\begin{Pro}
\label{pbfTh} Let G be a polycyclic-by-finite group containing a poly-infinite cyclic subgroup $H$ of finite index $n$. If $n$
is invertible in a regular ring $R$, then the group ring $R G$ is
also regular.
\end{Pro}
\begin{proof}
The ring $R H$ is regular \cite[1.5.11, 7.7.5]{McC}. The preceding lemma implies that $R G$ is also regular.
\end{proof}
\begin{Th}\label{33}
Let G be a polycyclic-by-finite group containing a 
poly-infinite cyclic subgroup of finite index $n$. Let $\al$, $\beta$ and $\gamma$ be group automorphisms such that $\al$ is of finite order $m$ and $\beta \circ \gamma$ is of finite order $m'$.
\begin{enumerate}
\item The group $\Nil_i(\zz G;\al)$ is an $(n \cdot m)$-primary torsion
group for $i \in \zz$.
\item The group $\Nil_i(\zz G;\zz G_{\beta},\zz G_{\gamma})$ is an
$(n \cdot m')$-primary torsion group for $i \in \zz$.
\end{enumerate}
\end{Th}
\begin{proof}
In the following, $\Nil_i(\zz G; \al)$ is treated, almost the same
arguments work for $\Nil_i(\zz G;\zz G_{\beta \circ \gamma}, \zz G)$. To obtain the general statement use the identity \cite[Proposition 3.2]{CP}:
\[
\Nil_i(R G; R G_\beta , R G_\gamma) \cong \Nil_i(R G; R G_{\beta \circ \gamma} , R G).\]

Since $\al$ is of finite order $m$, we can form the semidirect
product of $G$ and the cyclic group with $m$ elements $C_m$.
Conjugation with the element $(\text{id},1)$ extends $\al$ to an
inner automorphism of $G \rtimes C_m$, which is denoted by
$\tilde{\al}$. Since $\tilde{\al}$ is an inner automorphism, we
can apply Theorem~\ref{localization}. We have
$$
\zz[1/n,1/m] \ot_{\zz} \Nil_i\big( \zz G \rtimes C_m;\tilde{\al}
\big) \cong \Nil_i \big( \zz[1/n,1/m] G \rtimes C_m ;\tilde{\al}
\big).
$$
The right hand side vanishes since $\zz[1/n,1/m] G \rtimes C_m$ is regular (Proposition \ref{pbfTh}). Thus $\Nil_i(\zz G \rtimes C_m;
\tilde{\al})$ is $(n \cdot m)$-torsion. 
Lemma \ref{Transfer} implies now the theorem.
\end{proof}
\begin{Th}
Let G be an arbitrary group and let $\al$, $\beta$ and $\gamma$ be group
automorphisms such that $\al$ is of finite order $m$ and $\beta \circ \gamma$ is of finite order $m'$.
\begin{enumerate}
\item The group $\zz[1/m] \ot_{\zz} \Nil_i(\qq G;\al)$ is a
$\qq$-module for $i \in \zz$. 
\item The group $\zz[1/m'] \ot_{\zz} \Nil_i(\qq G;\qq
G_{\beta},\qq G_{\gamma})$ is a $\qq$-module for $i \in \zz$. 
\end{enumerate}
\end{Th}
\begin{proof}
We prove the result for Farrell Nil-groups almost the same
arguments work for the Waldhausen Nil-groups. We use the notation of the proof of Theorem~\ref{33}.
Since $\qq$ is a $\lambda$-ring Theorem~\ref{localization} implies that $\Nil_i\big( \qq G \rtimes C_m;\tilde{\al}\big)$ is a $\qq$-module. Lemma \ref{Transfer} implies now that $\zz[1/m] \ot_{\zz} \Nil_i(\qq G;\al)$ is a direct summand of $\zz[1/m] \ot_{\zz} \Nil_i\big( \qq G \rtimes C_m;\tilde{\al}\big)$ and therefore also a $\qq$-module.
\end{proof}
\subsection{The Relative Assembly Map} \label{Relative assembly map}
In the final section, we prove that the relative assembly map from
the family of finite groups to the family of virtually cyclic
groups is rationally an isomorphism. The main ingredients of the
proof are the torsion results of the preceding section. Combined with the calculation of $H_i^G(E_{\mathcal{F}in}(G); \textbf{K}_\zz) \ot \qq$ \cite{Lue2} we obtain the corollary stated in the introduction.

Before we start discussing the effect of the torsion results on
the Farrell-Jones Conjecture, let us briefly recall the relevant
notions. A family of subgroups of a group $G$ is a collection of subgroups of $G$ that is closed under conjugation and finite intersections. We denote the family of finite cyclic subgroups by $\mathcal{FC}y$, the family of finite subgroups by $\mathcal{F}in$, the
family of virtually cyclic subgroup by $\mathcal{VC}yc$ and the
family of all subgroups by $\mathcal{A}ll$. For such a family $\mathcal{F}$ there is a classifying space $E_\mathcal{F}(G)$, it is characterized by the property that
for any $G$-CW-complex $X$, all whose isotropy groups belong to
$\mathcal{F}$, there is up to $G$-homotopy precisely one $G$-map
from $X$ to $E_{\mathcal{F}}(G)$. 

Suppose we are given a family of subgroups $\mathcal{F}$ and a
subfamily $\mathcal{F}' \subseteq \mathcal{F}$. By the universal
property of $E_{\mathcal{F}}(G)$ we obtain a map
$E_{\mathcal{F}'}(G) \ra E_{\mathcal{F}}(G)$, which is unique up to
$G$-homotopy. Thus for every $G$-homology theory
$\mathcal{H}_{\ast}^G$ we obtain a \emph{relative assembly map}
$$
A_{\mathcal{F}' \ra \mathcal{F}} \colon \mathcal{H}_{\ast}^G
(E_{\mathcal{F}'}(G)) \ra \mathcal{H}_{\ast}^G
(E_{\mathcal{F}}(G)).
$$

In \cite{DL} it is explained how the $K$-theory spectrum of a ring
$R$, in the sequel denoted by \textbf{K}$_R$, gives an equivariant
 homology theory, in the sequel denoted by $H_i^?(- ;
\textbf{K}_{R})$. The Farrell-Jones conjecture \cite{FJ1} predicts that the
\emph{assembly map}
$$
A_{\mathcal{VC}yc \ra \mathcal{A}ll} \colon
H_i^G(E_{\mathcal{VC}yc}(G); \textbf{K}_{R} ) \ra  H_i^G (
E_{\mathcal{A}ll}(G) ; \textbf{K}_R) \cong K_i(R G)
$$
is an isomorphism. Assuming the Farrell-Jones conjecture is true,
the computation of $K_i(RG)$ reduces to the computation of
$H_i^G(E_{\mathcal{VC}yc}(G); \textbf{K}_{R})$. For a survey on the Farrell-Jones conjecture see for example~\cite{LR}.
\begin{Th} \label{mainth}
Let $G$ be a group and let $i$ be a natural number. If $i \geq 0$
then the rationalized relative assembly map
$$
A_{\mathcal{F}in \ra \mathcal{VC}yc} \colon
H_i^G(E_{\mathcal{F}in}(G); \textbf{K}_\zz) \ot \qq \ra
H_i^G(E_{\mathcal{VC}yc}(G); \textbf{K}_\zz) \ot \qq
$$
is an isomorphism. For $i < 0$, the relative assembly map is an
isomorphism even integrally.
\end{Th}
\begin{proof}
The proof follows closely a proof of the statement that the
relative assembly map is an isomorphism for a regular ring $R$ in
which the orders of all finite subgroups of $G$ are invertible
\cite[Proposition 2.14]{LR}. Because of the Transitivity Principle
\cite[Theorem 2.9]{LR} we need to prove that the
$\mathcal{F}in$-assembly map is an isomorphism for virtually
cyclic groups $V$. As mentioned in the introduction we can assume
that either $V \cong H \rtimes \zz$ or $V \cong G_1 \ast_H G_2$
with finite groups $H$, $G_1$ and $G_2$. In both cases we obtain
long exact sequences involving the algebraic $K$-theory of the
constituents, the algebraic $K$-theory of $V$ and additional
Nil-groups \cite{FH, Gr1, Wa1, Wa2}. If $V$ is a virtually cyclic group of the first type
or a virtually cyclic group of the second type the Nil-groups
vanish rationally by Theorem \ref{33}. Thus we get a long exact sequences
\begin{center}
 \makebox[0pt]{
$ \xym{ \cdots \ar[r] & K_i(\zz H) \ot \qq \ar[r] & K_i(\zz H) \ot \qq
\ar[r] & }$ }
\end{center}
\begin{center}
 \makebox[0pt]{
 $\xym{ \ar[r] & K_i(\zz V) \ot \qq \ar[r] & K_{i-1}(\zz H)
\ot \qq \ar[r] & K_{i-1}(\zz H) \ot \qq \ar[r] & \cdots } $}
\end{center}
and
\begin{center}
\makebox[0pt]{$
\xym@C=2em{ \cdots \ar[r] & K_i(\zz H) \ot \qq \ar[r] & \big( K_i(\zz G_1)
\op K_i(\zz G_2) \big) \ot \qq \ar[r] & K_i(\zz V) \ot \qq \ar[r] &}
$}
\end{center}
\[\xym{ \ar[r] & K_{i-1}(\zz H) \ot \qq \ar[r] &
\big( K_{i-1}(\zz G_1) \op K_{i-1}(\zz G_2) \big) \ot \qq \ar[r] &
\cdots }.
\]
One obtains analogous exact sequences for the sources of the
various assembly maps from the fact that the sources are
equivariant homology theories and specific models for
$E_{\mathcal{F}in}(V)$. These sequences are
compatible with the assembly maps. The assembly maps for finite
groups $H$, $G_1$ and $G_2$ are bijective. Now a Five-Lemma
argument shows that also the one for $V$ is bijective.

We obtain the stronger statement for $i < 0$, because the lower
Nil-groups are known to vanish \cite{FJ2}.
\end{proof}
\end{document}